\pdfoutput=1
\documentclass[twocolumn,a4paper]{svjour3}          % twocolumn
\usepackage{color}
\definecolor{dblue}{rgb}{0,0,0.5}
\definecolor{dgreen}{rgb}{0,0.5,0}

\usepackage{graphicx}
\usepackage{amsmath, amsfonts,graphicx}
 \usepackage{mathptmx}      % use Times fonts if available on your TeX system
\usepackage{subfig}
\usepackage{paralist,multirow}
\usepackage{isomath,epstopdf}
\usepackage{mathtools}
\usepackage{float,hyperref}
\usepackage{marvosym} %for \Letter
\usepackage[scaled=0.9]{helvet}
\hypersetup{
    colorlinks=true,
    linkcolor=dblue,
    filecolor=magenta,      
    urlcolor=dgreen,
    citecolor=dblue
  }
\usepackage
[
        a4paper,% other options: a3paper, a5paper, etc
        left=1.5cm,
        right=2cm,
]
{geometry}
\usepackage{microtype} 
\renewcommand{\d}{\mathrm{d}}
\newcommand{\intphi}{\int\limits_{\varphi}}
\newcommand{\intsphi}{\int_{\varphi}}
\begin{document}
\allowdisplaybreaks
\tolerance 1414
%\emergencystretch
\title{Optimization along Families of Periodic and Quasiperiodic Orbits in Dynamical Systems with Delay}


\author{Zaid Ahsan       \and 
       Harry Dankowicz \and Jan Sieber %etc.
}


\institute{ Z. Ahsan (\Letter) \and H. Dankowicz   \at
              Department of Mechanical Science and Engineering,\\
							University of Illinois Urbana-Champaign, Urbana, U.S.A \\
              \email{zaid2@illinois.edu}           %  \\
              \and
							J. Sieber \at
							College of Engineering, Mathematics and Physical 
				     Sciences,\\University of Exeter,
						  Exeter, United Kingdom
}


\maketitle

\begin{abstract}
This paper generalizes a previously-conceived, continuation-based optimization technique for scalar objective functions on constraint manifolds to cases of periodic and quasiperiodic solutions of delay-differential equations. A Lagrange formalism is used to construct adjoint conditions that are linear and homogenous in the unknown Lagrange multipliers. As a consequence, it is shown how critical points on the constraint manifold can be found through several stages of continuation along a sequence of connected one-dimensional manifolds of solutions to increasing subsets of the necessary optimality conditions. Due to the presence of delayed and advanced arguments in the original and adjoint differential equations, care must be taken to determine the degree of smoothness of the Lagrange multipliers with respect to time.  Such considerations naturally lead to a formulation in terms of multi-segment boundary-value problems (BVPs), including the possibility that the number of segments may change, or that their order may permute, during continuation. The methodology is illustrated using the software package \textsc{coco} on periodic orbits of both linear and nonlinear delay-differential equations, keeping in mind that closed-form solutions are not typically available even in the linear case. Finally, we demonstrate optimization on a family of quasiperiodic invariant tori in an example unfolding of a Hopf bifurcation with delay and parametric forcing. The quasiperiodic case is a further original contribution to the literature on optimization constrained by partial differential BVPs.

\end{abstract}

\section{Introduction}
\label{intro}

The optimization of time-delay systems has been the subject of intensive research for many years. Such systems arise naturally in control applications where unmodeled actuator dynamics results in delays between input signals and actuator responses~\cite{sieber2008tracking}, car following models that account for driver reaction times~\cite{orosz2005bifurcations}, and machine tool dynamics due to the regenerative effect~\cite{tlusty1963stability}. The wide range of applications has motivated the development of novel techniques for their optimization. For example, G\"{o}llmann \emph{et al.}~\cite{gollmann2009optimal} used a formulation based on the Pontryagin minimum principle to derive necessary optimality conditions for optimal control problems with delays in state and control variables. The obtained equations were discretized and transformed into a large-scale nonlinear programming model, which was then solved using off-the-shelf solvers. In another investigation, Yusoff and Sims~\cite{yusoff2011optimisation} combined the semi-discretization method~\cite{insperger2002semi} for time-periodic delay equations with differential evolution to optimize a variable helix/pitch tool geometry for regenerative chatter mitigation. Their results were also validated experimentally, confirming the predicted significant improvements in chatter stability. This problem of optimal selection of parameters for subtractive manufacturing was also reported in~\cite{iglesias2019optimum,iglesias2014optimisation,wojciechowski2017optimisation}.  Liao \emph{et al.}~\cite{liao2014nonlinear} developed an optimization technique for periodic solutions of delay differential equations using the harmonic balance method and continuation techniques. They posed an amplitude optimization problem subject to the algebraic constraints obtained by substitution of a truncated Fourier representation in the governing equation along with the stability conditions. The sensitivity expressions were analytically derived, and the optimization problem was then solved for the unknown Fourier coefficients and the unknown parameters. The delayed Duffing oscillator was used to validate the methodology. 

%In~\cite{lei2013optimal}, optimal sampled-data vibration control for uncertain nonlinear delay systems is considered. 

 %In another investigation by Liao et al.~\cite{liao2014nonlinear}, optimization along with continuation is used to optimize the periodic solutions of the classical delayed duffing oscillator. Their technique is based on the harmonic balance method, where the solution is represented by a truncated Fourier series. The sensitivity of the Fourier coefficients with respect to the parameters of interest is analytically derived, and the optimization problem is solved for the unknown Fourier coefficients and the unknown parameters.

The calculus of variations serves as a useful tool for constrained optimization problems. Here, a Lagrangian functional is constructed by combining the objective function with the imposed constraints using Lagrange multipliers (adjoint variables) as coefficients. The vanishing of the variations of the Lagrangian with respect to the design variables and the Lagrange multipliers then yields the necessary optimality conditions for a stationary point. In general, these equations cannot be solved directly. Instead, nonlinear solvers may be applied to various finite-dimensional discretizations. A major challenge with this approach is the selection of a good initial guess which converges to the desired solution. A resolution built on principles of parameter continuation was originally proposed in the work of Kern{\'e}vez and Doedel~\cite{kernevez1987optimization}. There, a sequence of properly initialized stages of continuation along one-dimensional manifolds of solutions to a subset of the necessary optimality conditions was used to connect the local extremum to an initial solution guess with vanishing Lagrange multipliers. This methodology was recently revisited by Li and Dankowicz~\cite{li2017staged} and there cast in terms of partial Lagrangians relevant to the general context of constrained optimization of integro-differential boundary-value problems \emph{without delay}. Importantly, this work showed how the Lagrangian structure was consistent with a staged construction paradigm implemented in the software package \textsc{coco}. 

In this work, we generalize the successive continuation approach of Kern{\'e}vez and Doedel to optimization along families of periodic and quasiperiodic orbits in dynamical systems \emph{with delay}. We derive the necessary optimality conditions from a suitably constructed Lagrangian without first discretizing the governing equations and unknowns. This approach is in contrast to other studies~\cite{rubino2018adjoint}, in which the discretization of the governing equations is first carried out and then the Lagrangian is constructed based on the discretized equations. In our formulation, the Lagrange multipliers satisfy coupled, piecewise-defined, boundary-value problems with both delayed and advanced arguments. Depending on the imposed constraints, the Lagrange multipliers may be discontinuous or nonsmooth at the interval boundary points, naturally resulting in a multi-segment problem \cite{calver2017numerical}. 

We first motivate our interest and approach with the problem of optimization of the response amplitude of a harmonically-forced, scalar, linear, delay-differential equation in Sect.~\ref{sec:motivation}. The general framework for problems with single delays is then considered in Sect.~\ref{sec:framework}, first for periodic orbits and subsequently for families of two-dimensional quasiperiodic invariant tori. As discussed in detail, the latter optimization problem falls into the category of constrained optimization for partial differential equations (PDEs)~\cite{matthias,hinze2008optimization,olikara2016computation}, for which the necessary optimality conditions take the form of coupled, piecewise-defined PDEs with non-local coupling, as well as associated boundary and interval conditions representing periodicity in one dimension and rotation in the other. Subsections of Sect.~\ref{sec:framework} consider example applications to the search for a saddle of the response amplitude of a harmonically-forced Duffing oscillator subject to delayed feedback control and a geometric fold along a family of quasiperiodic trajectories for constant rotation number. Analysis using the \textsc{coco} software package validates the successive continuation approach, as well as the simultaneous discretization of the dynamic constraints and adjoint equations. A number of additional considerations and opportunities for future work are considered in the concluding section.

 \section{Motivating Example}\label{sec:motivation}
\begin{figure}[h]
\centering
\subfloat{\includegraphics[width=0.95\columnwidth,trim={0 {0.0\textwidth} 0 0.0\textwidth},clip]{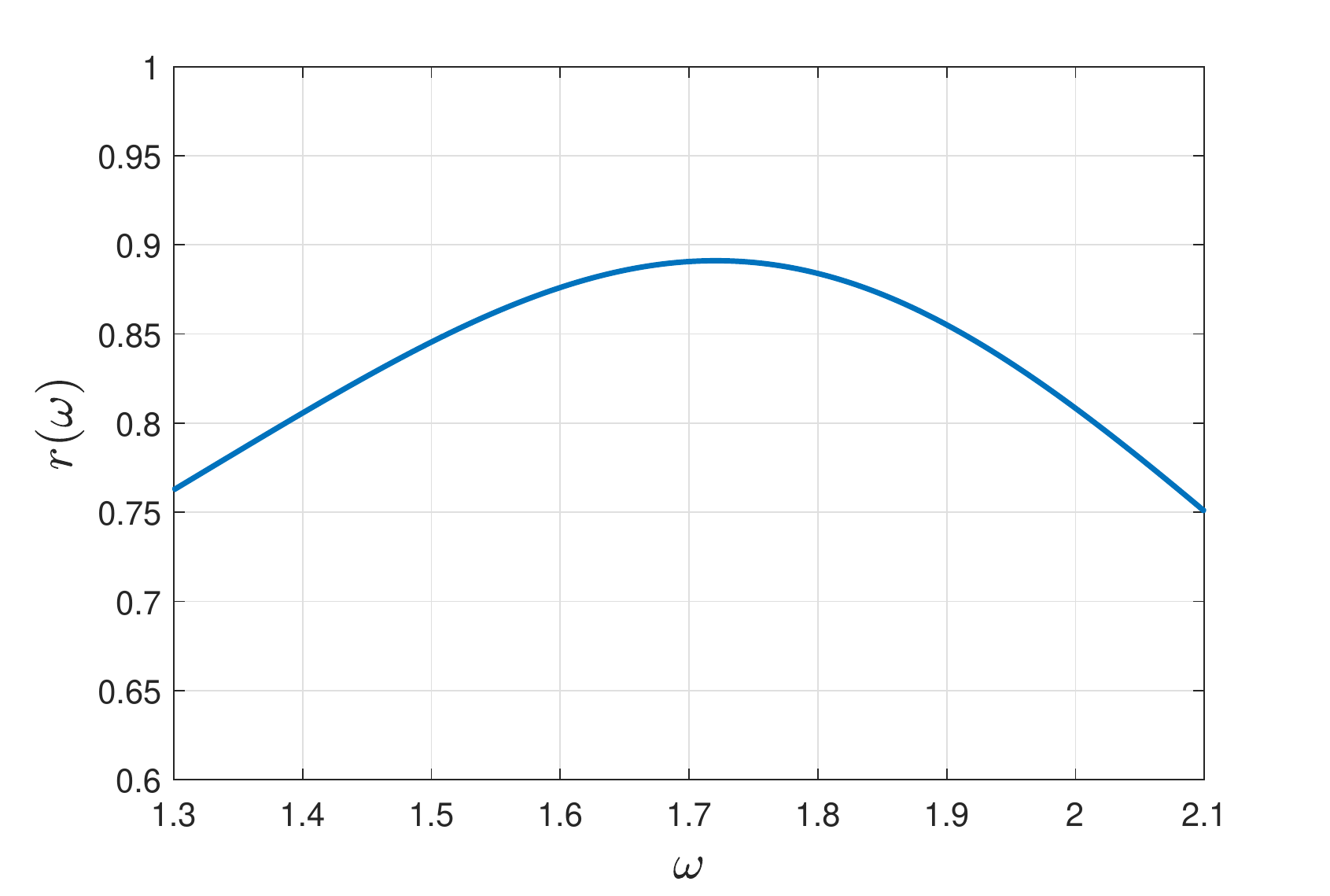}}
\caption{Frequency-response diagram for the steady-state periodic solutions of the harmonically-forced, scalar, linear delay-differential equation (\ref{ex1_eq1}). The maximum value of the amplitude is $r_\mathrm{crit}\approx 0.8911$ which occurs for $\omega=\omega_\mathrm{crit}\approx1.7207$ ($T=T_\mathrm{crit}\approx 3.6516$).}
\label{fig1}
\end{figure}

We illustrate the general framework for optimization along families of solutions to delay-differential equations (DDEs) by first considering periodic responses $z(t)$ of frequency $\omega$ for a harmonically-forced, scalar, linear, delay-differential equation
\begin{equation}
\label{ex1_eq1}
\dot{z}=-z- z\left(t-1\right)+ \text{cos}\omega t ,
\end{equation} 
where we omit (here, and throughout the paper) functional arguments when they are obvious from the context. It follows from the method of undetermined coefficients that such responses are of the harmonic form
\begin{equation}
z(t)=r(\omega)\cos(\omega t-\theta(\omega)),
\end{equation}
where
\begin{equation}
\label{ex_amp}
r(\omega)=\left[2+\omega ^2-2 \omega  \sin \omega +2 \cos \omega \right]^{-1/2}
\end{equation}
and
\begin{equation}
\cos\theta(\omega)=-\frac{1+\cos\omega}{r^3(\omega)},\,\sin\theta(\omega)=\frac{\sin\omega-\omega}{r^3(\omega)}.
\end{equation}
Let us consider the  optimization problem of finding the forcing frequency $\omega$ for which such a periodic response has maximum amplitude. It follows from \eqref{ex_amp} that the maximum amplitude $r_\mathrm{crit}\doteq r(\omega_\mathrm{crit})\approx 0.89$ is achieved for $\omega_{\mathrm{crit}}\approx 1.72$
(cf.~Fig.~\ref{fig1}), and that $z(t_\mathrm{crit})=r_\mathrm{crit}$ at time $t_{\mathrm{crit}}\doteq\theta(\omega_\mathrm{crit})/\omega_\mathrm{crit}\approx2.24$ (up to
  multiples of the period $T_\mathrm{crit}\doteq2\pi/\omega_\mathrm{crit}\approx3.65$). Hence, for this simple optimization problem all components of the solution are known exactly, enabling a comparison with the results of numerical algorithms.

\subsection{Formulation as a constrained optimization problem}
\label{sec:ex1_constrained}
We transform the above optimization problem into a  format suitable for a general numerical solver by introducing the excitation period $T=2\pi/\omega$ as an unknown ($T$ replaces $\omega$) and rescaling time (calling the new time $\tau$) such that $x(\tau)\doteq z(T\tau+T\phi/2\pi)$. Here, the free phase $\phi$ is to be chosen so as to shift the time on the interval $[0,1]$ when the periodic solution $x$ has a critical point to $\tau=0$. Thus, we are seeking a solution to the constrained optimization problem
\begin{equation}
  \label{ex1_mua}
  \mbox{maximize\quad}\mu_A=x(0)
\end{equation}
with respect to a continuous function $x$ on $[0,1]$, as well as the variables $T$ and $\phi$, subject to the equality constraints
\begin{align}
  x'&=-Tx-Tx\left(\tau+1-1/T\right)+T\cos\left(2 \pi \tau +\phi\right)\label{ex1_eq5a}\\
  &\mbox{\quad for $\tau\in(0,1/T)$,}\nonumber\\
  x'&=-Tx-Tx\left(\tau-1/T\right)+T\cos\left(2 \pi \tau +\phi\right)\label{ex1_eq5b}\\
  &\mbox{\quad for $\tau\in(1/T,1)$,}\nonumber\\
%\begin{subequations}
  0&=x\left(0\right)-x\left(1\right)\mbox{,}\label{ex1_eq5c}\\
  0&=x\left(0\right)+x\left(1-1/T\right)-\cos\phi. \label{ex1_eq5d}
%\end{subequations}
\end{align}
Here, the constraints \eqref{ex1_eq5a} and \eqref{ex1_eq5b} impose the original delay-differential equation on the interval $(0,1)$. They rely on periodicity to wrap the delayed argument back into this interval assuming that $T>1$. The constraints \eqref{ex1_eq5c} and \eqref{ex1_eq5d} are boundary conditions. Constraint \eqref{ex1_eq5c} imposes periodicity also on the interval boundary, while \eqref{ex1_eq5d} is a phase condition that ensures that $x'(0)=0$, consistent with $x$ having a critical point at $\tau=0$ and justifying the maximization of $x(0)$ as a substitute for the amplitude. By continuity of $x$ on $[0,1]$ and \eqref{ex1_eq5c} it follows that $x$ is, in fact, a smooth function on $[0,1]$. Indeed, from the explicit solution in the previous section, it follows that $x(\tau)=r(2\pi/T)\cos2\pi\tau$ and $\phi=\theta(2\pi/T)$ and, in particular, that optimality is obtained for $x(\tau)=x_\mathrm{crit}(\tau)\doteq r_\mathrm{crit}\cos2\pi\tau$ and $\phi=\phi_\mathrm{crit}\doteq\theta(2\pi/T_\mathrm{crit})$ for $T=T_\mathrm{crit}$.

The constrained optimization problem~\eqref{ex1_mua}--\eqref{ex1_eq5d} gives rise to the Lagrangian
\begin{align}
\label{ex1_eq9}
\lefteqn{L\left(x(\cdot),\phi,T, \mu_A, \lambda_1(\cdot),\lambda_2,\lambda_3,\eta_A\right)=\mu_A+ \eta_A\left(x\left(0\right)-\mu_A\right)}\nonumber\\
&+\int\limits_{0}^{1/T}\lambda_{1}\left[x^{\prime}+T\left[x+x\left(\tau+1-1/T\right)
-\cos\left(2\pi \tau+\phi\right)\right]\right]\,\mathrm{d}\tau\nonumber\\
&+\int\limits_{1/T}^{1-2/T}\lambda_{1}\left[x^{\prime}+T\left[x+x\left(\tau-1/T\right)
    -\cos\left(2\pi \tau+\phi\right)\right]\right]\,\mathrm{d}\tau\nonumber\\
&+\int\limits_{1-2/T}^{1-1/T}\lambda_{1}\left[x^{\prime}+T\left[x+x\left(\tau-1/T\right)
    -\cos\left(2\pi \tau+\phi\right)\right]\right]\,\mathrm{d}\tau\nonumber\\
&+\int\limits_{1-1/T}^{1}\lambda_{1}\left[x^{\prime}+T\left[x+x\left(\tau-1/T\right)
    -\cos\left(2\pi \tau+\phi\right)\right]\right]\,\mathrm{d}\tau\nonumber\\
&+\lambda_{2}\left(x\left(0\right)-x\left(1\right)\right)+\lambda_{3}\left(x\left(0\right)+x\left(1-1/T\right)-\text{cos}\left(\phi\right)\right),
\end{align} 
where the Lagrange multipliers are $\lambda_{1}(\tau)$ (a function on $[0,1]$) for the DDE constraints \eqref{ex1_eq5a} and \eqref{ex1_eq5b}, $\lambda_{2}$ and $\lambda_{3}$ for the boundary conditions \eqref{ex1_eq5c} and \eqref{ex1_eq5d}, and $\eta_A$ for the relationship between the fitness $\mu_A$ and $x(0)$ in \eqref{ex1_mua}. 

In \eqref{ex1_eq9}, the integral for the pairing between $\lambda_1$ and the DDE constraints has been split into $4$ parts, one for each of the intervals $(0,1/T)$, $(1/T,1-2/T)$, $(1-2/T,1-1/T)$, and $(1-1/T,1)$, reflecting different functional forms of the differential equations \eqref{ex1_eq5a} and \eqref{ex1_eq5b} for $x$ on $(0,1/T)$ and $(1/T,1)$, respectively, and anticipating possible discontinuities in $\lambda_1$ and $\lambda_1'$. For example, the split at $\tau=1-1/T$ is in anticipation of a potential discontinuity of the Lagrange multiplier $\lambda_1$ at this instant caused by the imposition of a constraint on $x$ evaluated at this time in \eqref{ex1_eq5d}. This discontinuity implies a potential discontinuity of $\lambda_1'$ at $\tau=1-2/T$. For the same reason, the appearances of $x(0)$ in \eqref{ex1_mua} and  \eqref{ex1_eq5d} suggest that $\lambda_1(0)\ne \lambda_1(1)$ resulting in a potential discontinuity of $\lambda_1'$ at $\tau=1-1/T$. All functions are assumed to be continuously differentiable on the partition implied by the integrals in \eqref{ex1_eq9}. The ordering of the discontinuity points assumes that $T>3$ (Fig.~\ref{fig1} shows that the optimal $T$ is in this range).

Imposing vanishing variations of the Lagrangian $L$ with respect to variations in all its arguments recovers the original constraints \eqref{ex1_mua}--\eqref{ex1_eq5d} and the following \emph{adjoint system} determining the Lagrange multipliers. Specfically, vanishing variations with respect to $x$ imply
\begin{equation}
\label{eq_adj1diffeq}
-\lambda_{1}^{\prime}+T\lambda_{1}+T\lambda_{1}\left(\tau+1/T\right)=0
\end{equation}
for $\tau\in\left(0,1/T\right)\cup\left(1/T,1-2/T\right)\cup\left(1-2/T,1-1/T\right)$ and
\begin{equation}
-\lambda_{1}^{\prime}+T\lambda_{1}+T\lambda_{1}\left(\tau-1+1/T\right)=0
\end{equation}
for $\tau\in\left(1-1/T,1\right)$. Boundary and interface conditions for these equations are obtained by considering variations with respect to $x(0)$, $x(1/T)$, $x(1-2/T)$, $x(1-1/T)$, and $x(1)$, corresponding in that order to
\begin{align}
0&=-\lambda_{1}(0)+\lambda_{2}+\lambda_{3}+\eta_A,\\
0&=\lambda_{1}(1/T)_--\lambda_{1}(1/T)_+,\label{eq_l1cont}\\
0&=\lambda_{1}(1-2/T)_--\lambda_{1}(1-2/T)_+,\label{eq_l1cont2}\\
0&=\lambda_{1}(1-1/T)_{-}-\lambda_{1}(1-1/T)_++\lambda_{3},\label{eq_l1discont}\\
0&=\lambda_{1}(1)-\lambda_{2},
\end{align}
using the convention that $\lambda_1(\tau^*)_\pm\doteq\lim_{{\tau\to\tau^*}\pm}\lambda_1(\tau)$ and recalling that $x(\tau)$ is continuous on $[0,1]$. Vanishing variations with respect to $\phi$ and $T$ imply the integral constraints
\begin{equation}
0=\int\limits_{0}^{1}T\lambda_{1}\sin\left(2\pi \tau+\phi\right)\,\mathrm{d}\tau+\lambda_{3}\sin\left(\phi\right)
\end{equation}
and
\begin{align}
\label{ex_last}
0&=\int\limits_0^{1/T}\lambda_1\left(x(\tau+1-1/T)+x'(\tau+1-1/T)/T\right)\,\mathrm{d}\tau\nonumber\\
&+\int\limits_{1/T}^1\lambda_1\left(x(\tau-1/T)+x'(\tau-1/T)/T\right)\,\mathrm{d}\tau\nonumber\\
&+\int\limits_0^1\lambda_1\left(x-\cos(2\pi\tau+\phi)\right)\,\mathrm{d}\tau+\lambda_3T^{-2}x'(1-1/T),
\end{align}
respectively. Finally, vanishing variation with respect to $\mu_A$ implies that
\begin{equation}
\label{algeqs}
1-\eta_A=0.
\end{equation}
In summary, the system of original contraints \eqref{ex1_mua}--\eqref{ex1_eq5d} and adjoint equations \eqref{eq_adj1diffeq}--\eqref{algeqs} is a nonlinear integro-differential boundary-value problem (BVP) defining the critical points of the Lagrangian $L$ and the constrained optimization problem \eqref{ex1_mua}--\eqref{ex1_eq5d}. 

In this example, the dimension of the manifold on which the constrained optimization problem is posed equals $1$, corresponding to the numbers of degrees of freedom of the nonlinear subsystem \eqref{ex1_mua}--\eqref{ex1_eq5d} (with variables $x$, $T$ and $\phi$). In contrast, the full system \eqref{ex1_mua}--\eqref{ex1_eq5d}, \eqref{eq_adj1diffeq}--\eqref{algeqs} has no such degrees of freedom and, consequently, generically has only isolated solutions. Several properties put it beyond the reach of ``off-the-shelf'' BVP solvers:
\begin{enumerate}
  \item\label{pt:ms} It consists of differential equations on multiple intervals (thus, the problem is called a \emph{multi-segment} BVP) with differential functional forms and continuous ``right-hand sides''. The number and length of these intervals is strongly problem  dependent, and may even change during the optimization process.
  \item The differential equations evaluate their right-hand sides at times deviating from $\tau$ (delayed or advanced arguments).
  \item\label{pt:nonlocal} The second point leads to nonlocal coupling across segments that is not restricted to coupling at the boundaries of the intervals. For example, \eqref{eq_adj1diffeq} couples the values of $\lambda_1$ in $(1-2/T,1-1/T)$ to values of $\lambda_1$ in $(1-1/T,1)$.
  \end{enumerate}

On the other hand,  the system \eqref{ex1_mua}--\eqref{ex1_eq5d}, \eqref{eq_adj1diffeq}--\eqref{algeqs} has some additional structure that aids both in its construction and solution:
\begin{enumerate}
  \item The equations are only forward coupled in that a solution to the original constraints \eqref{ex1_mua}--\eqref{ex1_eq5d} can be obtained independently of the values of the Lagrange multipliers, while a solution to the adjoint equations \eqref{eq_adj1diffeq}--\eqref{algeqs} requires knowledge of $x$, $T$, and $\phi$, and generically exists, at best, only for isolated choices of $x$, $T$, and $\phi$.
  \item The adjoint equations \eqref{eq_adj1diffeq}--\eqref{ex_last} (thus excluding \eqref{algeqs}) are linear and homogeneous in the Lagrange multipliers $\lambda_j$ ($j=1,2,3$) and $\eta_A$. A trivial solution of this subset of the adjoint system is therefore given by vanishing Lagrange multipliers for any $x$, $T$, and $\phi$.
    \item The adjoint equation \eqref{algeqs} is trivial both in construction and solution. Imposing its solution ($\eta_A=1$) on the remaining adjoint system, however, renders the latter nonhomogeneous.
\end{enumerate}
This structure will also be present for more general cases  than the example and can be exploited in the search for solutions, as well as to generate the adjoint equations \eqref{eq_adj1diffeq}--\eqref{ex_last} automatically during a staged construction of the optimization problem similar to \cite{li2017staged}.

In this example, a few facts about the Lagrange multipliers may be deduced directly from the adjoint equations. It follows immediately from \eqref{eq_l1cont} and \eqref{eq_l1cont2} that $\lambda_1$ is continuous at $\tau=1/T$ and $\tau=1-2/T$, and from \eqref{eq_adj1diffeq} that $\lambda_1'$ is continuous at $\tau=1/T$.  Moreover, using the explicitly known solution for $x$, it follows that the Lagrange multiplier $\lambda_3$ must equal $0$ at a local extremum. Indeed, substitution of the modified phase condition
\begin{equation}
  \delta=x\left(0\right)+x\left(1-1/T\right)-\cos\phi
\end{equation}
in lieu of \eqref{ex1_eq5d} implies that
\begin{equation}
\mu_A=\left(\cos\phi+\cos(\omega+\phi)+\omega\sin\phi\right)/r(\omega),
\end{equation}
where $\phi$ is implicitly determined by
\begin{equation}
\delta=\omega\left(\sin\phi+\sin(\omega+\phi)-\omega\cos\phi\right)/r(\omega)
\end{equation}
for $\delta\approx 0$. Implicit differentiation of both conditions with respect to the residual $\delta$ shows that the rate of change of $\mu_A$ with respect to $\delta$ equals $0$ at $\delta=0$. This, in turn, implies that that $\lambda_3=0$ at an extremum, i.e., that $\lambda_1$ is, in fact, continuous also at $\tau=1-1/T$ and, consequently, continuously differentiable also at $\tau=1-2/T$. In contrast, $\lambda'_1$ experiences a discontinuity at $\tau=1-1/T$ for nonzero $\eta_A=\lambda_1(0)-\lambda_1(1)$.

\subsection{Simple continuation}
\label{sec:ex1_success}

According to the properties enumerated above, a solution to \eqref{ex1_mua}--\eqref{ex1_eq5d}, \eqref{eq_adj1diffeq}--\eqref{algeqs} may be sought using a method of \emph{successive continuation} \cite{kernevez1987optimization,li2017staged} with an embedded multi-segment boundary-value problem implementation that permits evaluation of the right-hand side at arguments shifted by arbitrary times. Specifically, this method overcomes the problem of initializing a nonlinear solver for the full system by defining a sequence of continuation problems with one-dimensional solution manifolds that connect an initial solution guess with Lagrange multipliers all equal to $0$ with the sought critical point for which $\eta_A$ must equal $1$.

% Specifically, \js{successive continuation addresses}
% the problem of initializing a nonlinear solver by defining a sequence
% of continuation problems with one-dimensional solution manifolds that
% connect an initial solution guess with Lagrange multipliers all equal
% to $0$ with \js{a critical point of the optimization problem.}

To this end, we consider the system given by the relationship between $\mu_A$ and $x(0)$ in \eqref{ex1_mua}, the
boundary-value problem constraints \eqref{ex1_eq5a}--\eqref{ex1_eq5d}, and the adjoint integral-differential boundary-value problem \eqref{eq_adj1diffeq}--\eqref{ex_last}, but purposely omit the algebraic constraint \eqref{algeqs}. 
Although we anticipate that $\lambda_3$ will equal $0$ throughout the analysis, we keep $\lambda_3$ as an unknown and monitor its value during continuation. By linearity and homogeneity of the adjoint subsystem in the Lagrange multipliers $\lambda_j$ and $\eta_A$, it follows that solutions to the full system lie on either of two one-dimensional manifolds. The first of these consists of functions $x(\tau)=r(2\pi/T)\cos2\pi\tau$ with corresponding $T$, $\phi$, and $\mu_A=x(0)=r(2\pi/T)$, and with vanishing Lagrange multipliers. The second manifold consists of the periodic solution $x_\mathrm{crit}(\tau)=r_\mathrm{crit}\cos2\pi\tau$ with corresponding $T=T_\mathrm{crit}$, $\phi=\phi_\mathrm{crit}$, and $\mu_A=r_\mathrm{crit}$, and with varying Lagrange multipliers proportional to $\eta_A$. The two manifolds clearly intersect at the local extremum of $\mu_A$ along the first manifold. The sought solution to the complete set of equations \eqref{ex1_mua}--\eqref{ex1_eq5d}, \eqref{eq_adj1diffeq}--\eqref{algeqs} corresponds to the point along the second manifold where $\eta_A=1$.

In this example, the solutions along the first manifold are known explicitly. In other cases, an initial periodic response may be approximately obtained from the dynamically stable solution by direct simulation. Given such an initial solution guess for $x$, $T$, and $\phi$, a nonlinear solver may be employed to converge to a point on the manifold. A numerical continuation algorithm (e.g., pseudo-arclength continuation) may then be used to generate a sequence of points along the manifold, meanwhile monitoring for local extrema of $\mu_A$ and singular points for the system Jacobian (corresponding to branch points on the manifold). As shown above, and true also in the general case, these coincide. Using standard techniques, numerical continuation may proceed from such a branch point along the secondary manifold with the help of a candidate direction of continuation, for example, one that is i) transversal to the tangent direction to the original solution manifold and ii) in the plane spanned by the tangent directions to the two manifolds at the branch point.

\begin{figure}[h]
\centering
\subfloat{\includegraphics[width=0.9\columnwidth,trim={0 {0.0\textwidth} 0 0.0\textwidth},clip]{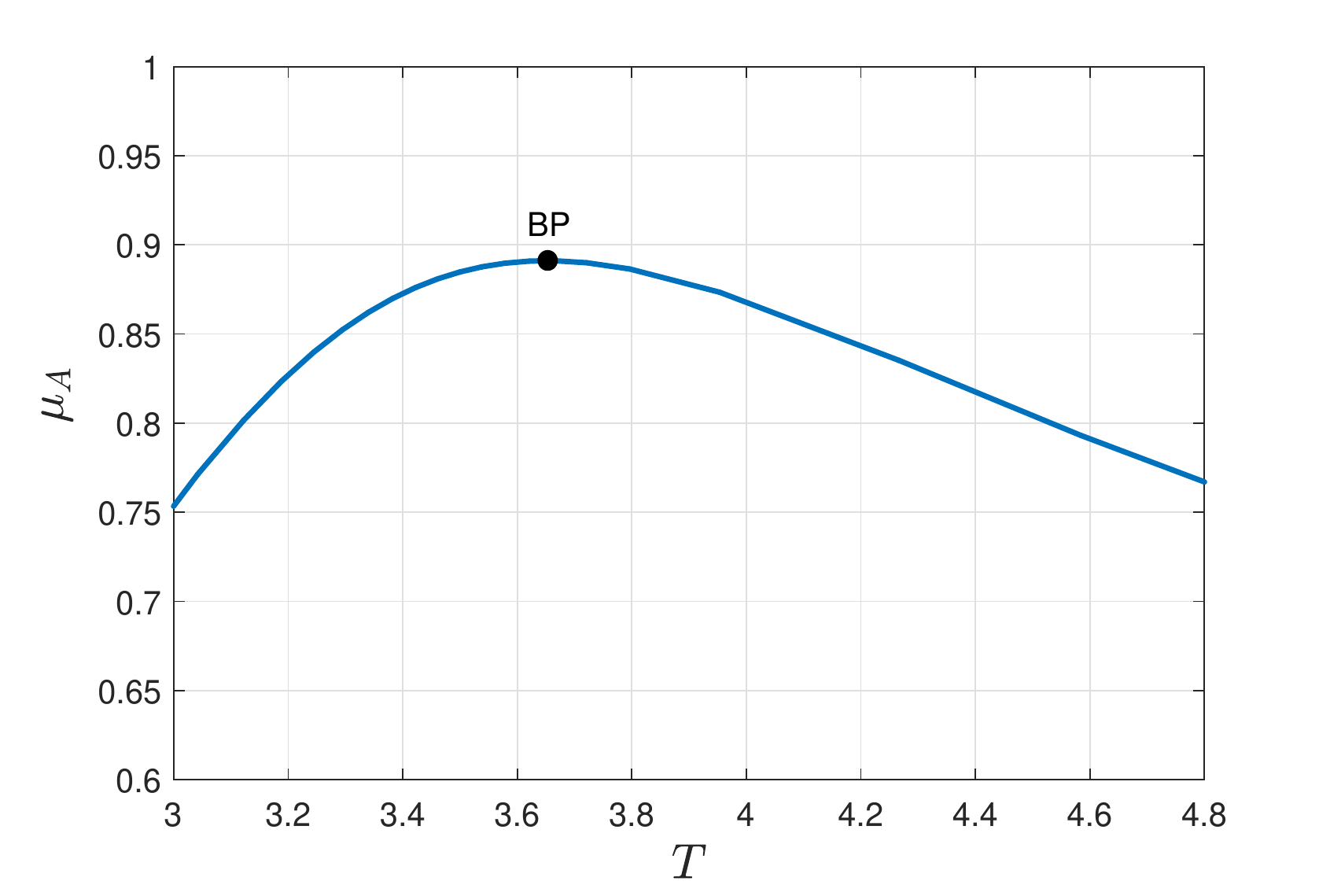}}
\caption{Results from numerical continuation with vanishing Lagrange multipliers. The maximum value of $\mu_A$ is located at $T\approx 3.6515$ and is here identified by the label BP, since it approximately coincides with a branch point.}
\label{fig2}
\end{figure}

\begin{figure}[h]
\centering
\subfloat[]{\includegraphics[width=0.9\columnwidth,trim={0 {0.0\textwidth} 0 0.0\textwidth},clip]{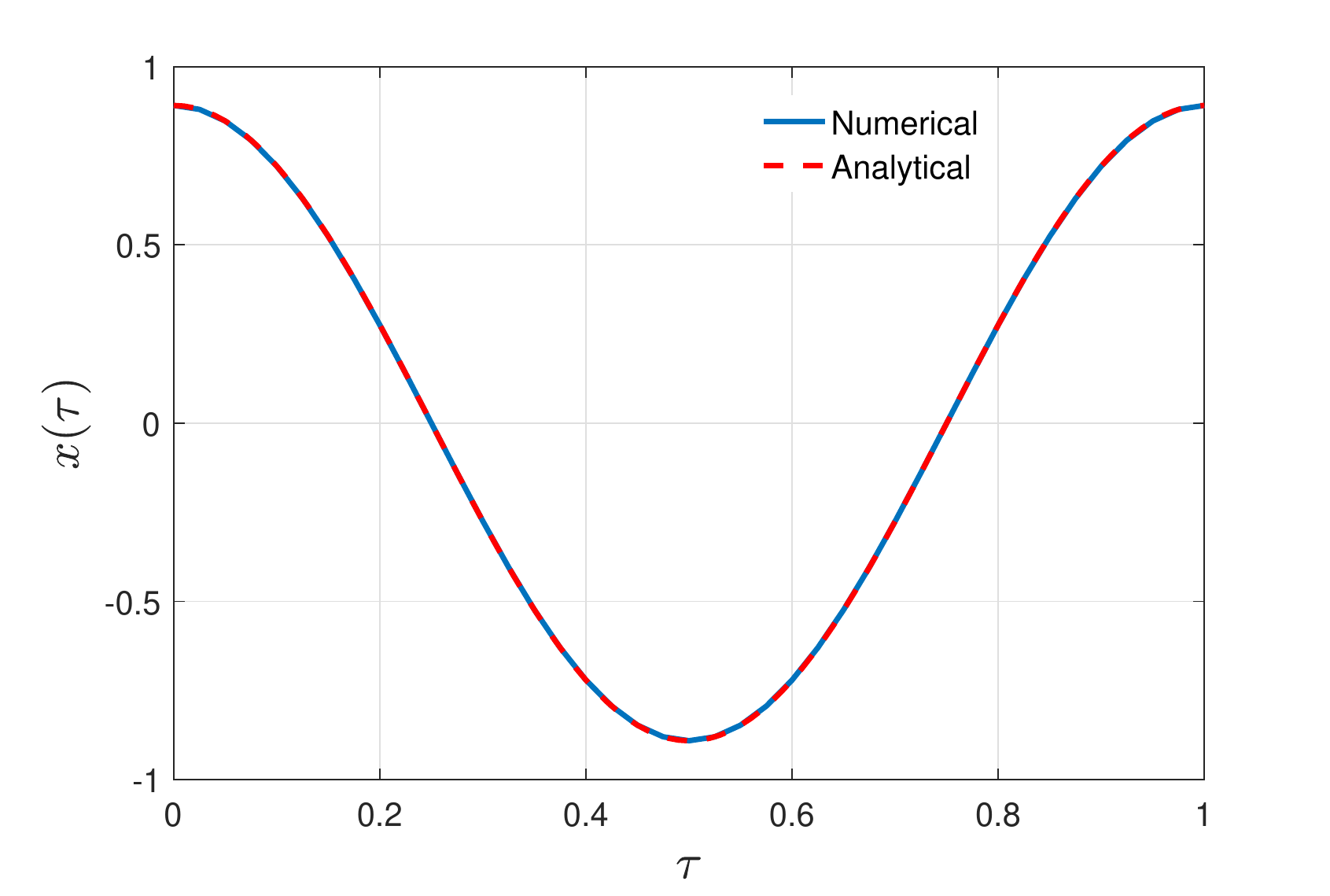}}\\
\subfloat[]{\includegraphics[width=0.9\columnwidth,trim={0 {0.0\textwidth} 0 0.0\textwidth},clip]{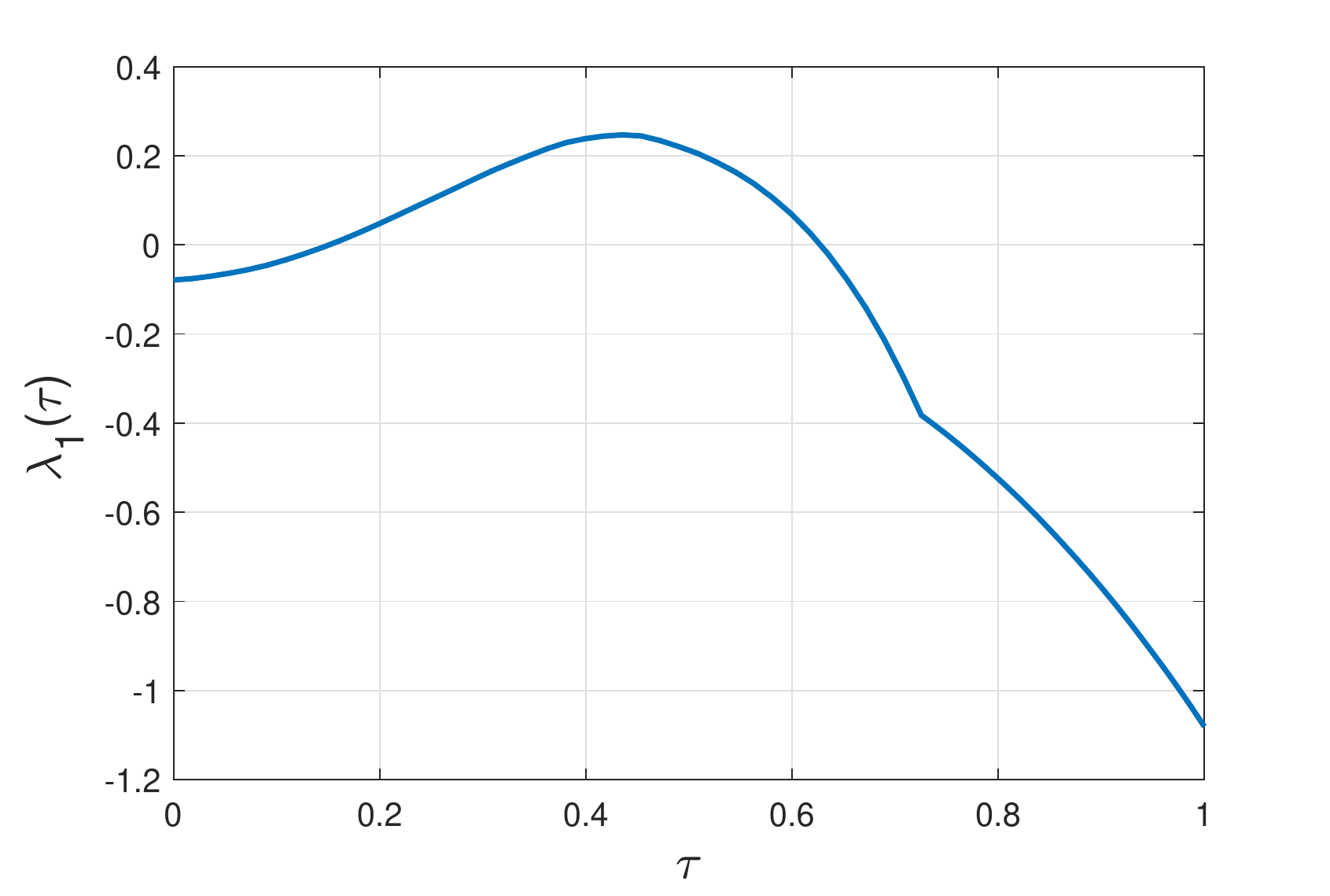}}
\caption{(a) $x(\tau)$ and (b) $\lambda_1(\tau)$ at the terminal point of the second stage of continuation with $\eta_A=1$. The upper panel shows a comparison between the numerical solution and the analytical solution at the extremum. The bottom panel shows the Lagrange multiplier associated with the imposition of the DDE admitting a slope discontinuity at $\tau=1-1/T$.}
\label{fig3}
\end{figure}

Continuation using such an implementation in the \textsc{coco} software package~\cite{COCO} approximately locates an extremum (in the form of a fold point in $\mu_A$ along the solution manifold) at $T\approx 3.6515$ as shown in Fig.~\ref{fig2}. Branch switching from the nearby branch point (exact coincidence is lost due to discretization) and continuation until $\eta_A=1$ yields the graphs of $x(\tau)$ and $\lambda_1(\tau)$ shown in Fig.~\ref{fig3}. As seen in the bottom panel, $\lambda_1(\tau)$ is approximately continuous at $1-1/T\approx 0.73$, albeit with discontinuous derivative at this point, since $\lambda_1(0)-\lambda_1(1)=1$.

\section{General Optimization Framework}\label{sec:framework}
In this section, we discuss the general methodology for optimization on periodic and quasiperiodic solutions $z(t)\in\mathbb{R}^{n}$ of delay-differential equations with a single delay of the form
\begin{equation}
\label{field}
\dot{z}=f\left(t,z,z\left(t-\alpha\right),p\right),
\end{equation}
where $f:\mathbb{R}^{1}\times\mathbb{R}^{n}\times\mathbb{R}^{n}\times\mathbb{R}^{q}\rightarrow\mathbb{R}^{n}$ is periodic in its first argument with period $T$. Here, $\alpha$ and $p$ denote the time delay and the problem parameters (excluding $T$), respectively. 

As the motivating example in the previous section illustrates, the problem Lagrangian and, by implication, the adjoint equations are linear in the Lagrange multipliers. The adjoint equations may therefore be constructed term-by-term by successively deriving the contributions from disjoint collections of constraints from the corresponding \emph{partial Lagrangians} associated with a subset of the Lagrange multipliers. Until the full set of constraints has been considered, the adjoint equations are not completely known. The following subsections discuss the partial Lagrangians and the implied contributions to the adjoint equations resulting from the DDE constraints and boundary conditions associated with periodic and quasiperiodic orbits. For particular examples, we indicate the additional contributions associated with problem-specific constraints that complete the construction of the adjoint equations.  In all cases, the contribution from the objective function to the Lagrangian implies the algebraic adjoint condition that the corresponding Lagrange multiplier ($\eta_A$ in the previous section) must equal $1$ at a stationary point.

\subsection{Periodic orbits}\label{formulation}
\label{sec_periodic}
Suppose first that $T>\alpha$ and consider the problem of optimizing a scalar-valued objective functional on a family of continuous solutions $x(\tau)$ to the differential equations
\begin{align}\label{per:dde1}
x^{\prime}&=Tf\left(T\tau,x,x\left(\tau+1-\alpha/T\right),p\right)&&\mbox{for $\tau \in \left(0,\alpha/T\right)$,}\\
\label{per:dde2}
x^{\prime}&=Tf\left(T\tau,x,x\left(\tau-\alpha/T\right),p\right)&&\mbox{for $\tau \in \left(\alpha/T,1\right)$,}
\end{align}
 and the boundary conditions
\begin{equation}
x\left(0\right)-x\left(1\right)=0.
\label{ms_eq2}
\end{equation} 
By a rescaling of the independent variable by $T$, such solutions correspond to periodic solutions of \eqref{field} with period $T$. By continuity and periodicity, such solutions must be continuously differentiable to all orders.

Suppose, in fact, that $T>3\alpha$ and that the objective functional and any additional constraints depend on pointwise values of $x(\tau)$ only at $\tau=0$, $\tau=1$, and $\tau=\beta$ for some $\beta=\beta(\alpha,T)$ such that
\begin{equation}
2\alpha/T<\beta<1-\alpha/T.
\end{equation} 
As we show below, such dependence results in an additional adjoint equation associated with variations with respect to $x(\beta)$. Other pointwise dependencies of the objective functional would be treated similarly, while dependence on an integral over the entire interval $[0,1]$ of a function of $x$ would not result in additional adjoint equations. We may formulate a corresponding partial Lagrangian
\begin{align}
&L_\mathrm{BVP}\left(x(\cdot),\alpha,T,p,\lambda_f(\cdot),\lambda_\mathrm{bc}\right)=\nonumber\\
&\quad\lambda_{\mathrm{bc}}^\text{T}\left(x(0)-x(1)\right)+\int\limits_{0}^{\alpha/T}\lambda_{f}^\text{T}\left(x^{\prime}-Tf_{1}\right)\,\mathrm{d}\tau \nonumber\\
 &\quad+\int\limits_{\alpha/T}^{\beta-\alpha/T}\lambda_{f}^\text{T}\left(x^{\prime}-Tf_{0}
 \right)\,\mathrm{d}\tau+\int\limits_{\beta-\alpha/T}^{\beta}\lambda_{f}^\text{T}\left(x^{\prime}-Tf_{0}
 \right)\,\mathrm{d}\tau \nonumber\\
&\quad+\int\limits_{\beta}^{1-\alpha/T}\lambda_{f}^\text{T}\left(x^{\prime}-Tf_{0} \right)\,\mathrm{d}\tau +\int\limits_{1-\alpha/T}^{1}\lambda_{f}^\text{T}\left(x^{\prime}-Tf_{0} \right)\,\mathrm{d}\tau,\label{eq_Lpartial}
\end{align}
where $f_{j}(\tau)=f\left(T\tau,x\left(\tau\right),x\left(\tau+j-\alpha/T\right),p\right)$. Here, $\lambda_{f}(\tau)$ and $\lambda_{\mathrm{bc}}$ are the Lagrange multipliers associated with the imposition of the differential equations and boundary conditions, respectively, and each integrand is assumed to be continuously differentiable on the corresponding interval. The splitting of the integral is here motivated by an anticipated discontinuity of $\lambda_f$ at $\tau=\beta$ and, consequently, of $\lambda_f'$ at $\tau=\beta-\alpha/T$, the different functional forms of the original DDEs on the intervals $(0,\alpha/T)$ and $(\alpha/T,1)$, and an anticipated discontinuity in $\lambda_f'$ also at $\tau=1-\alpha/T$.

By the stated assumptions on the objective function and any additional constraints, it is easy to show that, at a stationary point of the total Lagrangian, $\lambda_f(\tau)$ must be continuous at $\tau=\alpha/T$, $\tau=\beta-\alpha/T$, and $\tau=1-\alpha/T$. Using the notation
\begin{align}
f_{j,k}(\tau)&=\partial_kf(T\tau,x(\tau),x(\tau+j-\alpha/T),p),\\
f_{j,q}(\tau)&=\frac{\d}{\d\, q}f\left(T\tau,x\left(\tau\right),x\left(\tau+j-\alpha/T\right),p\right)
\end{align}
for $j=0,1$ and $q=\alpha,T$ ($\partial_kf$ is the partial derivative of $f$ with respect to its $k$th argument, $\d/\d q$ is the total derivative of an expression with respect to $q$), the contributions to the necessary adjoint conditions for a stationary point of the total Lagrangian are given by
\begin{align}
&-\lambda_{f}^{\prime\text{T}}-T\lambda_{f}^{\text{T}}f_{1,2}-T\lambda_{f}^{\text{T}}\left(\tau+\alpha/T\right)f_{0,3}\left(\tau+\alpha/T\right)
\label{per:adj:dde1}
\end{align}
for variations with respect to $x(\tau)$ on $\tau\in\left(0,\alpha/T\right)$;
\begin{align}
\label{eq_middleadjointDDE}
&-\lambda_{f}^{\prime\text{T}}-T\lambda_{f}^{\text{T}}f_{0,2}-T\lambda_{f}^{\text{T}}\left(\tau+\alpha/T\right)f_{0,3}\left(\tau+\alpha/T\right)
\end{align}
for variations with respect to $x(\tau)$ on $\tau\in\left(\alpha/T,\beta-\alpha/T\right)\cup\left(\beta-\alpha/T,\beta\right)\cup\left(\beta,1-\alpha/T\right)$;
\begin{align}
\label{eq_finaladjointDDE}
&-\lambda_{f}^{\prime\text{T}}-T\lambda_{f}^{\text{T}}f_{0,2}-T\lambda_{f}^{\text{T}}\left(\tau+\alpha/T-1\right)f_{1,3}\left(\tau+\alpha/T-1\right)
\end{align}
for variations with respect to $x(\tau)$ on $\tau\in\left(1-\alpha/T,1\right)$;
\begin{align}\label{per:adj:bd}
&-\lambda_{f}^\text{T}\left(0\right)+\lambda_\mathrm{bc}^\text{T},
&&\lambda_{f}^\text{T}(\beta)_--\lambda_{f}^\text{T}(\beta)_+,
&&\lambda_{f}^\text{T}\left(1\right)-\lambda_\mathrm{bc}^\text{T},
\end{align}
for variations with respect to $x(0)$, $x(\beta)$, and $x(1)$, respectively;
\begin{align}
&-\int_{0}^{\alpha/T}\lambda_{f}^{\text{T}}Tf_{1,\alpha}\,\mathrm{d}\tau-\int_{\alpha/T}^{1}\lambda_{f}^{\text{T}}Tf_{0,\alpha}\,\mathrm{d}\tau
\label{eq_varalpha}
\end{align}
for variations with respect to $\alpha$;
\begin{align}
&-\int_{0}^{\alpha/T}\lambda_{f}^{\text{T}}\left(f_{1} +Tf_{1,T}\right)\,\mathrm{d} \tau-\int_{\alpha/T}^{1}\lambda_{f}^{\text{T}}\left(f_{0}+Tf_{0,T}\right)\,\mathrm{d}\tau
\label{eq_varT}
\end{align}
for variations with respect to $T$; and
\begin{align}
&-\int_{0}^{\alpha/T}\lambda_{1}^{\text{T}}Tf_{1,4}\,\mathrm{d}\tau-\int_{\alpha/T}^{1}\lambda_{1}^{\text{T}}Tf_{0,4}\,\mathrm{d}\tau
\label{per:adj:p}
\end{align}
for variations with respect to $p$. The terms $f_{j,T}$ and $f_{j,\alpha}$ in \eqref{eq_varalpha} and \eqref{eq_varT} both contain time derivatives $x'$ with delayed or advanced arguments, since $T$ and $\alpha$ both appear in the evaluation of $x$ in the third arguments of $f_0$ and $f_1$.

As previously anticipated, the explicit dependence of the Lagrangian on the internal state point $x(\beta)$ results in a potential discontinuity of the Lagrange multiplier $\lambda_{f}(\tau)$  at $\tau=\beta$. Continuous differentiability of $x(\tau)$ on $[0,1]$ and of $\lambda_f(\tau)$ on $(0,\beta-\alpha/T)$, $(\beta-\alpha/T,\beta)$, $(\beta,1-\alpha/T)$, and $(1-\alpha/T,1)$ implies that the necessary conditions for an extremum are in the form of a multi-segment boundary-value problem in a single trajectory segment for $x(\tau)$ and four coupled trajectory segments for $\lambda_f(\tau)$. A similar result is obtained, for example, in the limiting case when $\beta=1-\alpha/T$. This case specializes to the example discussed in the previous section, since there $\alpha=1$, $\beta=1-1/T$, and $T>3$. In contrast, when $\beta$ is either $0$ or $1$, i.e., when there is no dependence of the objective function or any additional constraints on an internal point, then we obtain a single trajectory segment for $x(\tau)$ and three coupled trajectory segments for $\lambda_f(\tau)$ with both variables continuous throughout the interval $[0,1]$. 

\subsection{A Duffing oscillator with delayed PD control}
\label{sec:duffing}
As an application of the general methodology when $\beta=0$, consider the harmonically-forced Duffing oscillator with delayed state (proportional and derivative; PD) feedback given by the DDE
\begin{align}
&\ddot{z}+2\zeta\dot{z}+z+\mu z^{3}=2 a z\left(t-\alpha\right)\nonumber\\
&\qquad \qquad\qquad\qquad+2 b \dot{z}\left(t-\alpha\right)+\gamma\,\cos \left(2\pi t/T\right).
\end{align}
Inspired by \cite{hu1998resonances}, for fixed $\zeta$, $\mu$, $a$, $b$, and $\gamma$, we seek a delay $\alpha$ that minimizes the maximum amplitude of oscillation along a family of periodic responses of this system under variations in the excitation period $T$.  Since the optimization problem involves minimizing a maximum, it corresponds to the search for a saddle point in the value of the oscillation amplitude on the two-dimensional constraint manifold.

Following Section~\ref{sec:ex1_constrained}, let $x_1(\tau)\doteq z(T\tau+T\phi/2\pi)$ and $x_2(\tau)\doteq \dot{z}(T\tau+T\phi/2\pi)$ represent the displacement and velocity, respectively, on the rescaled time interval $[0,1]$. The phase $\phi$ is again to be chosen so as to shift the time on this interval when the oscillator reaches its maximum displacement to $\tau=0$. It follows that the objective functional is given by
\begin{equation}
\mu_A=x_1(0)
\label{per:amp}
\end{equation}
for solutions of \eqref{per:dde1}--\eqref{ms_eq2} subject to the phase condition 
\begin{equation}
x_2(0)=0
\label{per:phase}
\end{equation}
and corresponding to the vector field
\begin{align}
&f(t,u,v,p)=\begin{pmatrix}u_2\\-2\zeta u_2-u_1-\mu u_1^3\end{pmatrix}\nonumber\\
&\qquad\qquad\quad+\begin{pmatrix}0\\2av_1+2bv_2+\gamma\cos\left(2\pi t/T+\phi\right)\end{pmatrix},
\end{align}
where $p=\phi$.

The partial  Lagrangian for the objective functional and phase condition is
\begin{align}
L_\mathrm{opt}\left(x(\cdot),\mu_A,\lambda_\mathrm{ph},\eta_A\right)=
\mu_A+\lambda_\mathrm{ph}x_2(0)+\eta_A(x_1(0)-\mu_A),
\end{align}
where $\lambda_\mathrm{ph}$ and $\eta_A$ are additional Lagrange multipliers. This partial Lagrangian adds the term $(\eta_A, \lambda_\mathrm{ph})^\text{T}$ to the variation with
  respect to $x(0)$ in \eqref{per:adj:bd} (first term) and results in the algebraic adjoint constraint 
  \begin{align}\label{per:adj:mua}
    0&=1-\eta_A,
  \end{align}
assuming no additional dependence of the problem Lagrangian on $\mu_A$.

Since neither the objective functional nor the additional phase condition depend on $x$ evaluated at an interior point of the interval $[0,1]$, it follows that $\lambda_f$ is continuous on the entire interval. This simplifies the partial Lagrangian $L_\mathrm{BVP}$ in \eqref{eq_Lpartial} as the two integrals with boundary $\beta$ can be combined, and the resulting adjoint DDE contribution \eqref{eq_middleadjointDDE} can be applied on the combined interval $(\alpha/T,1-\alpha/T)$, such that $\lambda_f(\tau)$ is in fact continuously differentiable on $(0,1-\alpha/T)$. Correspondingly, the middle adjoint condition in \eqref{per:adj:bd} can be omitted. Moreover, like in the motivating example in Section~\ref{sec:ex1_constrained}, it is easy to see that the rate of change of $\mu_A=x_1(0)$ with respect to $\delta=x_2(0)$ vanishes at $\delta=0$. We conclude that $\lambda_\mathrm{ph}=0$ at a stationary point of the Lagrangian. This implies that $\lambda_{f,2}(1)=\lambda_{f,2}(0)$ and, by inspection of \eqref{eq_middleadjointDDE} and \eqref{eq_finaladjointDDE}, that both components of $\lambda_f$ are actually continuously differentiable throughout the interval $[0,1]$.

Since the dimension $n_\mathrm{opt}$ of the optimization manifold equals $2$, the successive continuation approach proposed by Kern{\'e}vez and Doedel~\cite{kernevez1987optimization} requires multiple stages (in contrast to the motivating example in Section~\ref{sec:ex1_constrained}, where $n_\mathrm{opt}=1$): one initially optimizes only with respect to one variable, following a curve in the optimization manifold, keeping $n_\mathrm{opt}-1$ variables fixed. At each successive stage of continuation one releases one further optimization variable, until all variables are free. In this analysis, we propose to keep $\alpha$ fixed during the initial stage of continuation, corresponding to the imposition of a constraint on the set of unknowns. To this end, we consider the additional partial Lagrangian
\begin{equation}
  \label{per:lcont}
  L_\mathrm{sc}(\alpha,\mu_\alpha,\eta_\alpha)=\eta_\alpha(\alpha-\mu_\alpha)\mbox{.}
\end{equation}
This partial Lagrangian adds the constraint
\begin{equation}
  \label{per:mualpha}
  \alpha=\mu_\alpha
\end{equation}
and the algebraic adjoint equation (for vanishing variation with respect to $\mu_\alpha$)
\begin{align}\label{per:adj:mualpha}
  0&=\eta_\alpha,
\end{align}
and adds $\eta_\alpha$ to the adjoint variations with respect to $\alpha$ in \eqref{eq_varalpha}. The total problem Lagrangian is now given by 
\begin{align}
&L\left(x(\cdot),\alpha,T,p,\mu_A,\mu_\alpha, \lambda_f(\cdot),\lambda_\mathrm{bc},\lambda_\mathrm{ph},\eta_A,\eta_\alpha\right)=\nonumber\\
&\qquad L_\mathrm{BVP}\left(x(\cdot),\alpha,T,p,\lambda_f(\cdot),\lambda_\mathrm{bc}\right)\nonumber\\
&\qquad\qquad +L_\mathrm{opt}\left(x(\cdot),\mu_A,\lambda_\mathrm{ph},\eta_A\right)+L_\mathrm{sc}(\alpha,\mu_\alpha,\eta_\alpha).
\end{align}
The necessary conditions for an extremum of the total Lagrangian are then given by (i) the original differential equations and boundary
conditions, \eqref{per:dde1}--\eqref{ms_eq2}, \eqref{per:amp}, \eqref{per:phase}, and \eqref{per:mualpha}, and (ii) the various adjoint equations, including \eqref{per:adj:mua} and \eqref{per:adj:mualpha}, assembled in stages as constraints and variables are added, setting the sums of all resulting contributions equal to $0$. Although we anticipate that $\lambda_\mathrm{ph}$ will equal $0$ throughout the analysis, we keep $\lambda_\mathrm{ph}$ as an unknown and monitor its value during continuation.

As in the previous section, we may locate an extremum of $L$ by several successive stages of continuation, in each stage omitting one or both of the adjoint conditions \eqref{per:adj:mua} and \eqref{per:adj:mualpha}. In particular, by holding $\mu_\alpha$ fixed and letting $\eta_A$ vary freely, we may arrive at a solution with $\eta_A=1$ in two stages of continuation: first, by continuing along a one-dimensional manifold with vanishing Lagrange multipliers, and next by branch-switching at a local extremum of $\mu_A$ to a secondary branch along which only the Lagrange multipliers vary and, in fact, do so proportionally to $\eta_A$. A final stage of continuation then proceeds from the point on this second manifold where $\eta_A=1$, but this time with $\eta_A$ fixed at $1$ and $\mu_\alpha$ free to vary. A sought extremum is obtained when $\eta_\alpha=0$.

%Li \emph{et al.}\ \cite{li2017staged} make this gradual release possible in the continuation package \textsc{COCO} by introducing continuation parameters for each optimization variable that is initially inactive, adding one further partial Lagrangian for each initially inactive optimization variable. In our case, this is (\emph{inactive} in the terminology of \cite{dankowicz2013recipes,li2017staged}). 

\begin{figure}[h]
\centering
\subfloat{\includegraphics[width=0.90\columnwidth,trim={0 {0.0\textwidth} 0 0.0\textwidth},clip]{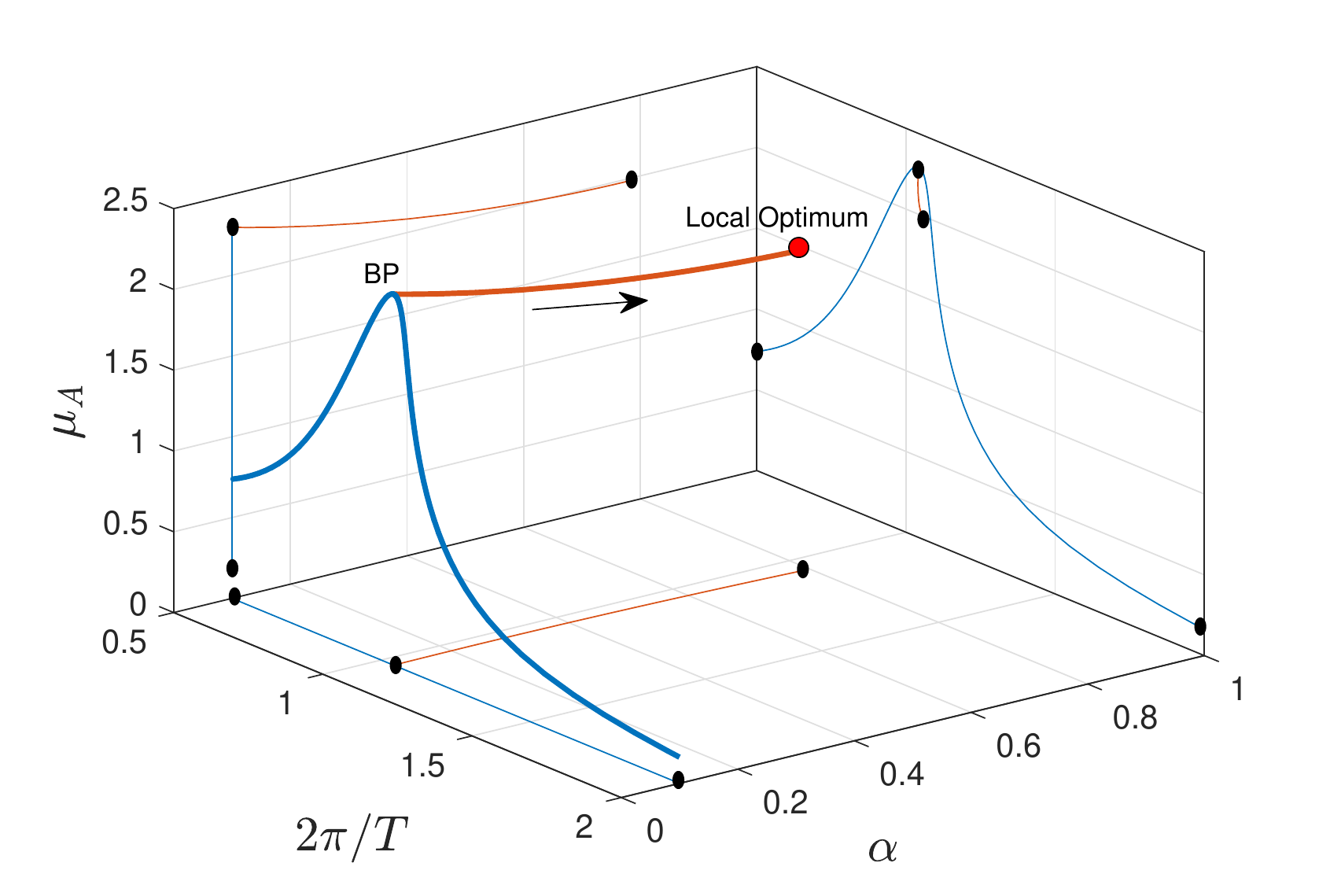}}
\caption{Optimization of the displacement amplitude along periodic orbits of the harmonically-excited Duffing oscillator with $\zeta=0.05$, $\mu=0.05$, $a=0.05$, $b=-0.05$, and $\gamma=0.5$ under variations in $\alpha$ and $T$. Three successive stages of continuation connect the sought saddle point with an initial solution guess with vanishing Lagrange multipliers. Stages 1 (blue) and 3 (red) described in the text are visible in the  $(\alpha,2\pi/T,\mu_A)$ space. In Stage 1, a peak in the displacement amplitude is approximately detected in close proximity to a branch point for the corresponding continuation problem. The second stage involves branch switching to a branch along which only the Lagrange multipliers vary (not visible). The red curve shows the final stage of continuation with fixed $\eta_{A}=1$. The optimal delay and corresponding period obtained at the terminal point with $\eta_\alpha=0$ equal $0.7824$ and $5.88$, respectively. At this point $\mu_A=1.9852$.}
\label{fig4}
\end{figure}

An example of such an analysis for the case when $\zeta=0.05$, $\mu=0.05$, $a=0.05$, $b=-0.05$, and $\gamma=0.5$ is shown in Fig.~\ref{fig4} (projected into the $(\alpha, 2\pi/T,\mu_A)$ space). Here, the full integro-differential boundary-value problem is discretized and analyzed using the \textsc{coco}~\cite{COCO} package following the methodology discussed in~\cite{dankowicz2013recipes} in terms of continuous, piecewise-polynomial approximants on a uniform partition of every solution segment into $N=10$ mesh intervals, resulting in a large system of nonlinear algebraic equations. The successive continuation approach then proceeds along the following stages:
\begin{itemize}
\item \emph{Initial guess.} An initial solution guess for $x(\tau)$ near the first manifold is first constructed using direct simulation with $\alpha=0.1$ and $T=2\pi$, after which $\phi$ is adjusted such that the maximum of $x_1$ occurs at $\tau\approx0$. We finally let $\mu_A=x_1(0)$ and $\mu_\alpha=\alpha$.
\item \emph{Stage 1: Continuation along manifold with vanishing Lagrange multipliers.} The delay $\alpha$ is held constant by fixing $\mu_\alpha$ at its initial value. Continuation proceeds along the blue curve in Fig.~\ref{fig4}, monitoring for branch points (coincident with extrema in $\mu_A$ up to discretization errors).
\item\emph{Stage 2: Continuation along manifold with varying Lagrange multipliers.} Branch off at the discovered branch point (labeled BP in Fig.~\ref{fig4}) with $\mu_\alpha$ still fixed, stopping when $\eta_A$ reaches $1$. During this continuation all primary variables $x(\cdot)$, $\phi$, $T$, $\alpha$ stay constant. Only Lagrange multipliers change their values. This continuation does not change any coordinates in Fig.~\ref{fig4} (we remain at the point BP).
\item\emph{Stage 3: Continuation with varying $\mu_\alpha$.} Fix $\eta_A$ at $1$ and allow  $\mu_\alpha$ (and, consequently, $\alpha$) to vary. Continue while monitoring $\eta_\alpha$ for zero crossings (red curve in Fig.~\ref{fig4}). The point where $\eta_\alpha=0$ along the red curve is labelled ``Local Optimum''. At this point all necessary conditions for a stationary point of $L$ are satisfied, including $\eta_A=1$ and $\eta_\alpha=0$.
\end{itemize}

\begin{figure}[h]
\centering
\subfloat[]{\includegraphics[width=0.9\columnwidth,trim={0 {0.0\textwidth} 0 0.0\textwidth},clip]{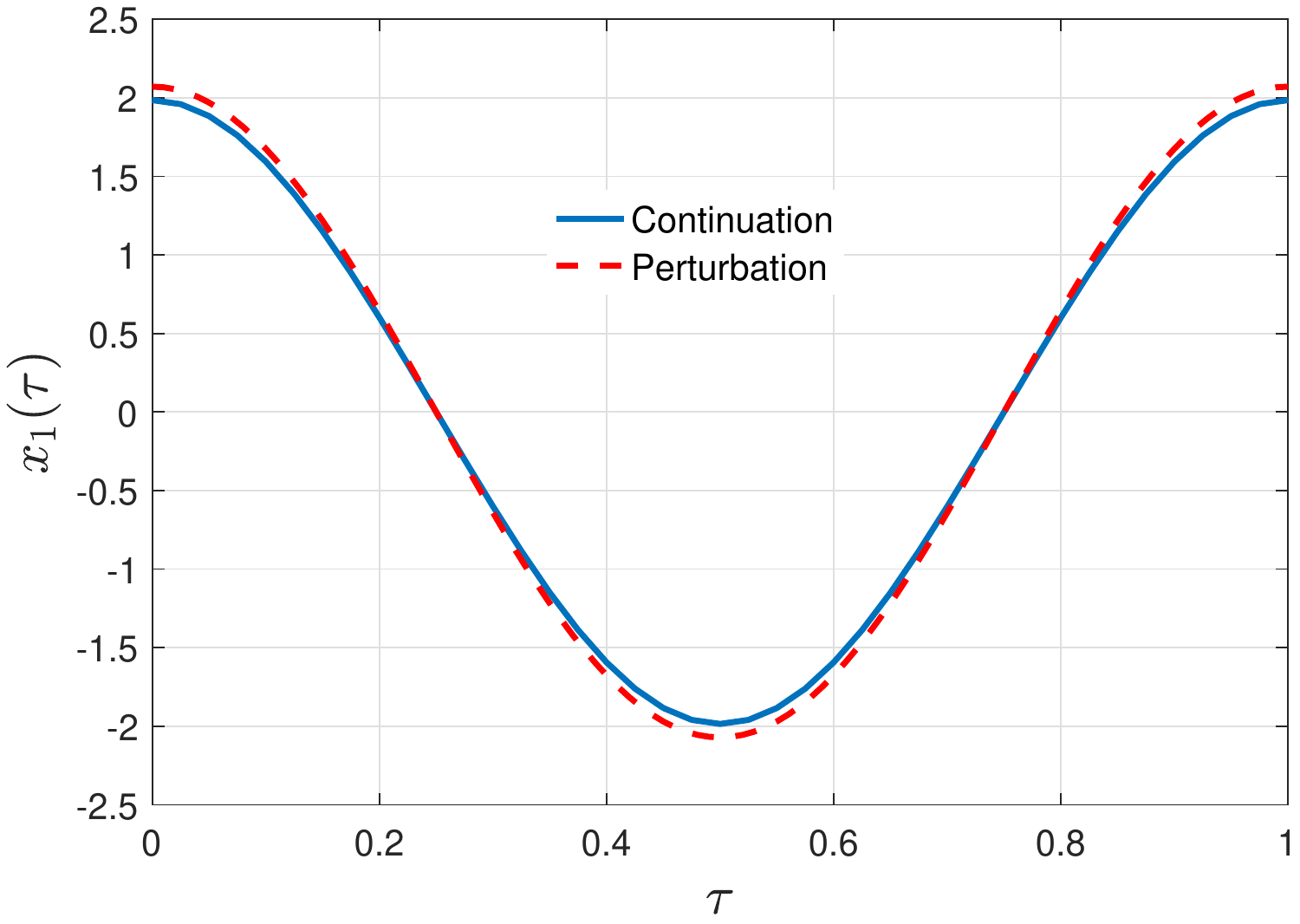}}\\
\subfloat[]{\includegraphics[width=0.9\columnwidth,trim={0 {0.0\textwidth} 0 0.0\textwidth},clip]{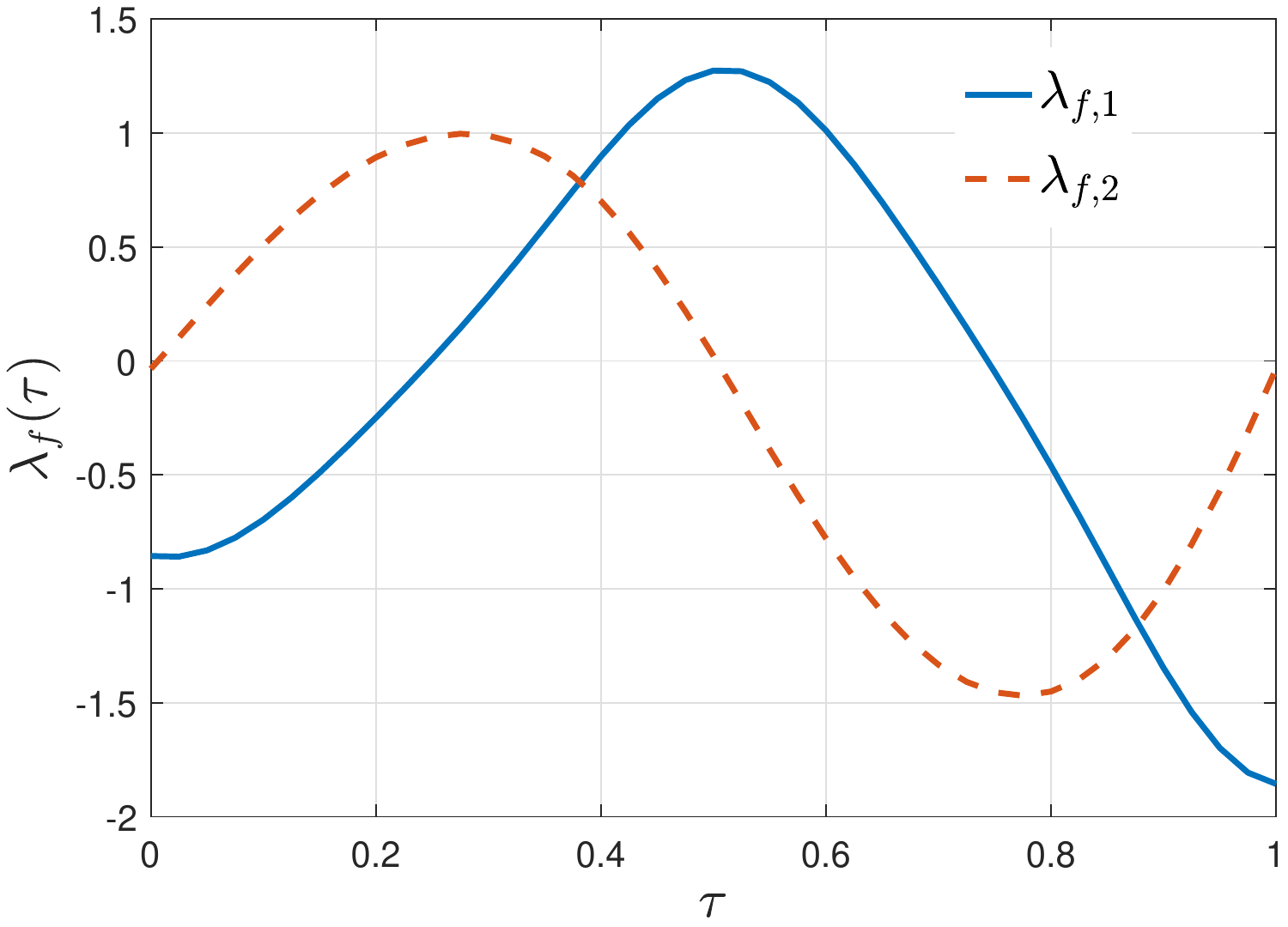}}
\caption{(a) $x_1(\tau)$ and (b) $\lambda_{f,1}(\tau)$ and $\lambda_{f,2}(\tau)$ at the terminal point of the third stage of continuation illustrated in Fig.~\ref{fig4}. The upper panel shows a comparison between the numerical solutions obtained using continuation at the computed optimal value of $\alpha$, with a first-order multiple-scales perturbation analysis at the predicted optimal value of $\alpha$.} 
\label{fig5}
\end{figure}

The end point of stage 3 corresponds to a critical point at $\alpha\approx0.7824$, $\phi\approx1.488$, and $T\approx 5.88$ (which Fig.~\ref{fig4} confirms to be a saddle point).  We may compare the resulting optimal delay with the prediction from a first-order multiple-scales perturbation analysis for the weakly nonlinear (small $\mu$), weakly damped (small $\zeta$), and weakly forced (small $\gamma$) case, which predicts a maximal (with respect to $T$) response amplitude
\begin{equation}\label{ms:amplitude}
  \cfrac{\gamma}{2|\zeta+a\sin\alpha-b\cos\alpha|}\mbox{,}
\end{equation}
(independent of $\mu$, see the appendix for intermediate steps and \cite{hu1998resonances,rusinek2014dynamics}). The computed optimal delay $\alpha\approx 0.7824$ is in close agreement with the predicted optimal delay $\pi/4\approx 0.7854$ obtained from \eqref{ms:amplitude} for the case that $b=-a$. The optimal displacement profile $x_1(\tau)$ and the components of $\lambda_f(\tau)$ are shown in Fig.~\ref{fig5}. The top panel shows close agreement between the results obtained using continuation and the harmonic response obtained from the perturbation analysis, at the computed optimal values of $\alpha$, $T$, and $\phi$. Panel (b) of Fig.~\ref{fig5} shows the functional Lagrange multipliers $\lambda_f$, confirming that they are approximately smooth in this example (since the objective does not depend on $\beta\in(0,1)$) but with $\lambda_{f,1}(1)\ne\lambda_{f,1}(0)$ and $\lambda_{f,2}(1)\approx\lambda_{f,2}(0)$ (since the objective functional and the phase constraint depend on $x(0)$ and $\lambda_\mathrm{ph}\approx 0$).

\begin{figure}
\centering
%\subfloat{\includegraphics[width=0.9\columnwidth,trim={0 {0.0\textwidth} 0 0.0\textwidth},clip]{fig7a}}\\
\subfloat[]{\includegraphics[width=0.9\columnwidth,trim={0 {0.0\textwidth} 0 0.0\textwidth},clip]{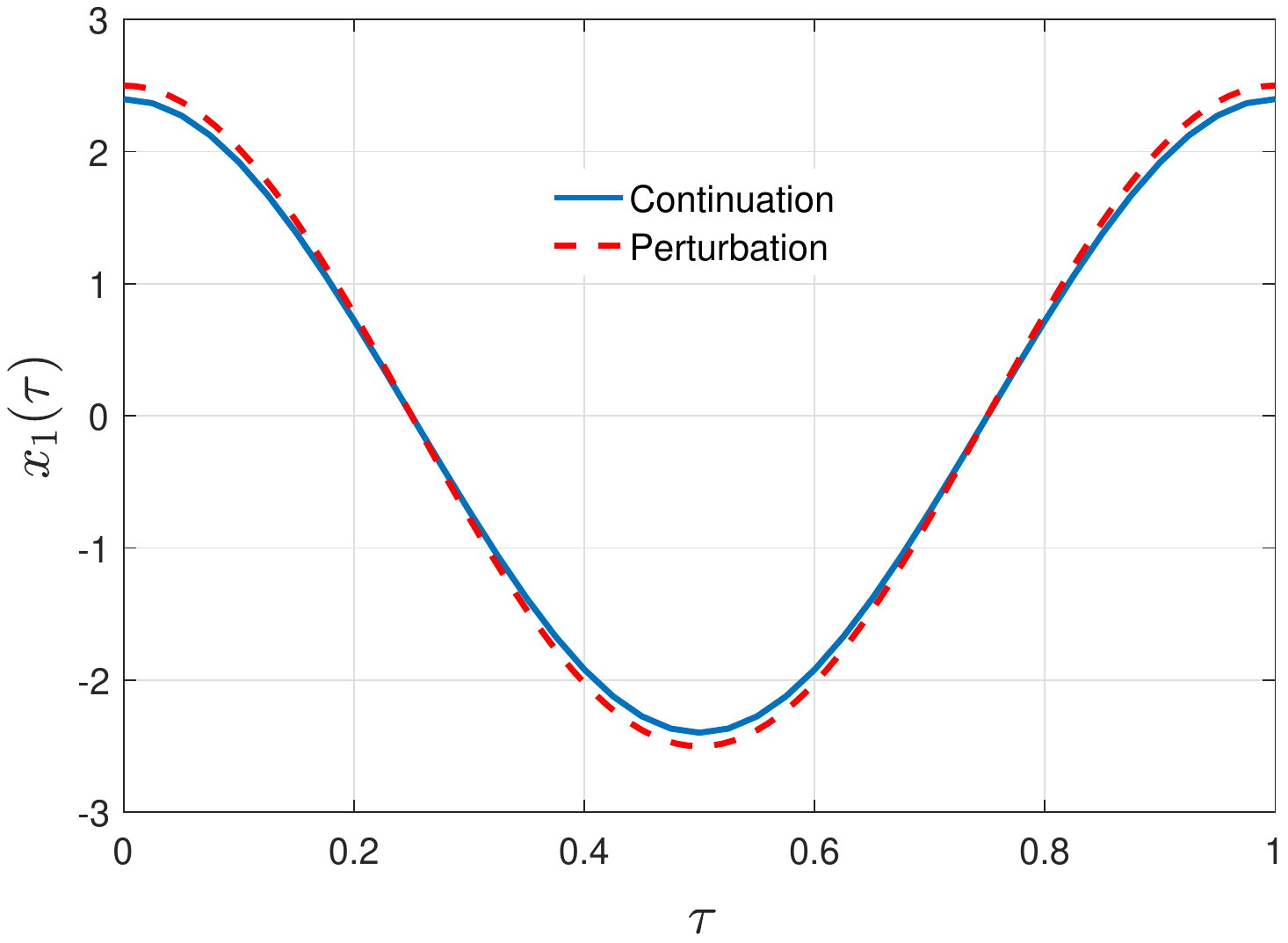}}\\
\subfloat[]{\includegraphics[width=0.9\columnwidth,trim={0 {0.0\textwidth} 0 0.0\textwidth},clip]{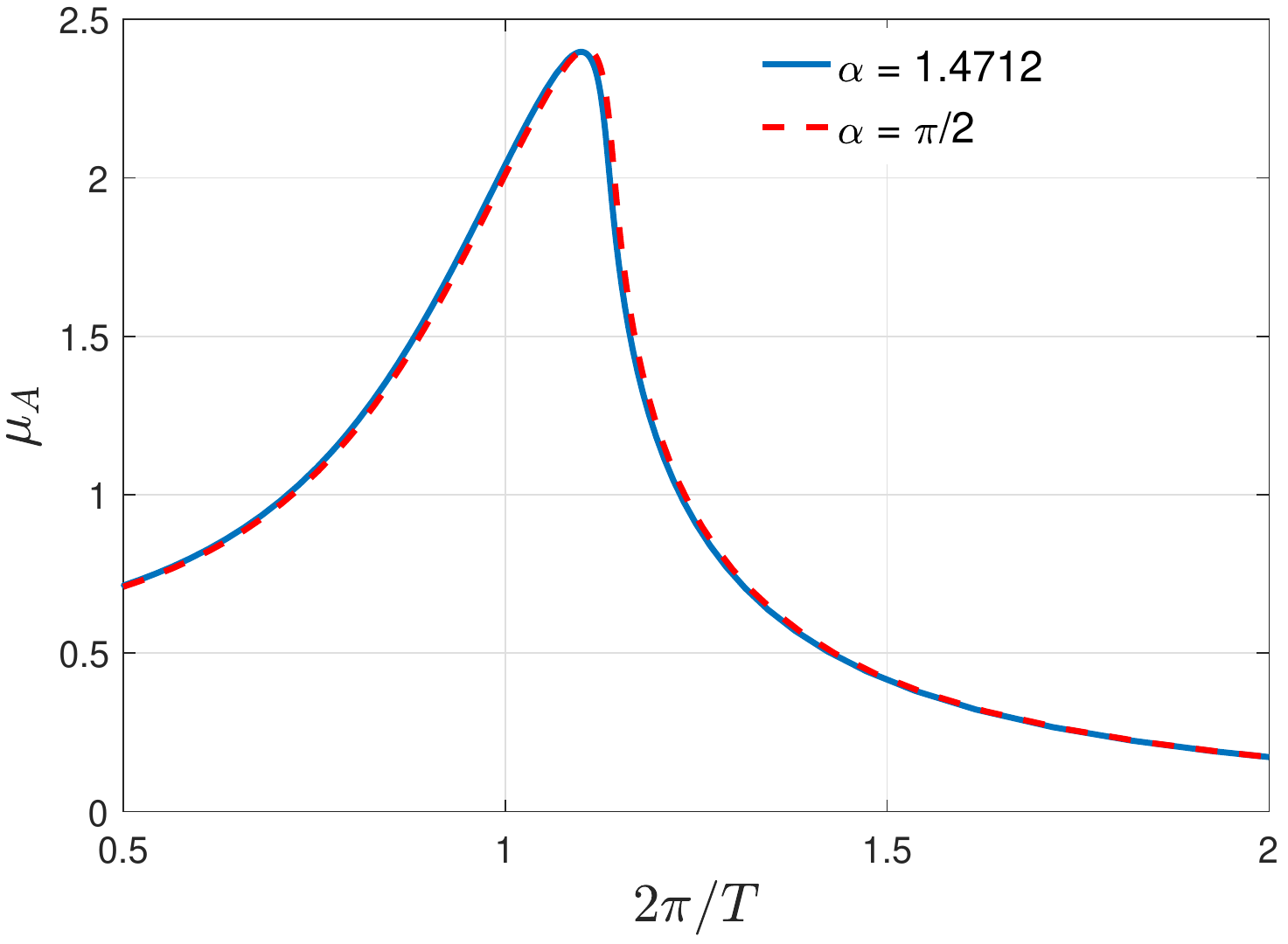}}
%\subfloat{\includegraphics[width=0.49\textwidth,trim={0 {0.0\textwidth} 0 0.0\textwidth},clip]{fig_ex2_4c}}
\caption{Optimization of the displacement amplitude along periodic orbits of the harmonically-excited, weakly-nonlinear Duffing oscillator with $\zeta=0.05$, $\mu=0.05$, $a=0.05$, $b=0$, and $\gamma=0.5$ under variations in $\alpha$ and $T$. (a) Comparison of the displacement profile obtained from continuation at the computed optimal delay $\alpha\approx1.4712$ and period $T\approx5.7151$ with the results predicted by perturbation analysis. (b) Frequency-response diagrams for the computed and predicted critical delay values $1.4712$ and $\pi/2$, respectively.}
\label{fig6}
\end{figure}

\begin{figure}
\centering
\subfloat[]{\includegraphics[width=0.9\columnwidth,trim={0 {0.0\textwidth} 0 0.0\textwidth},clip]{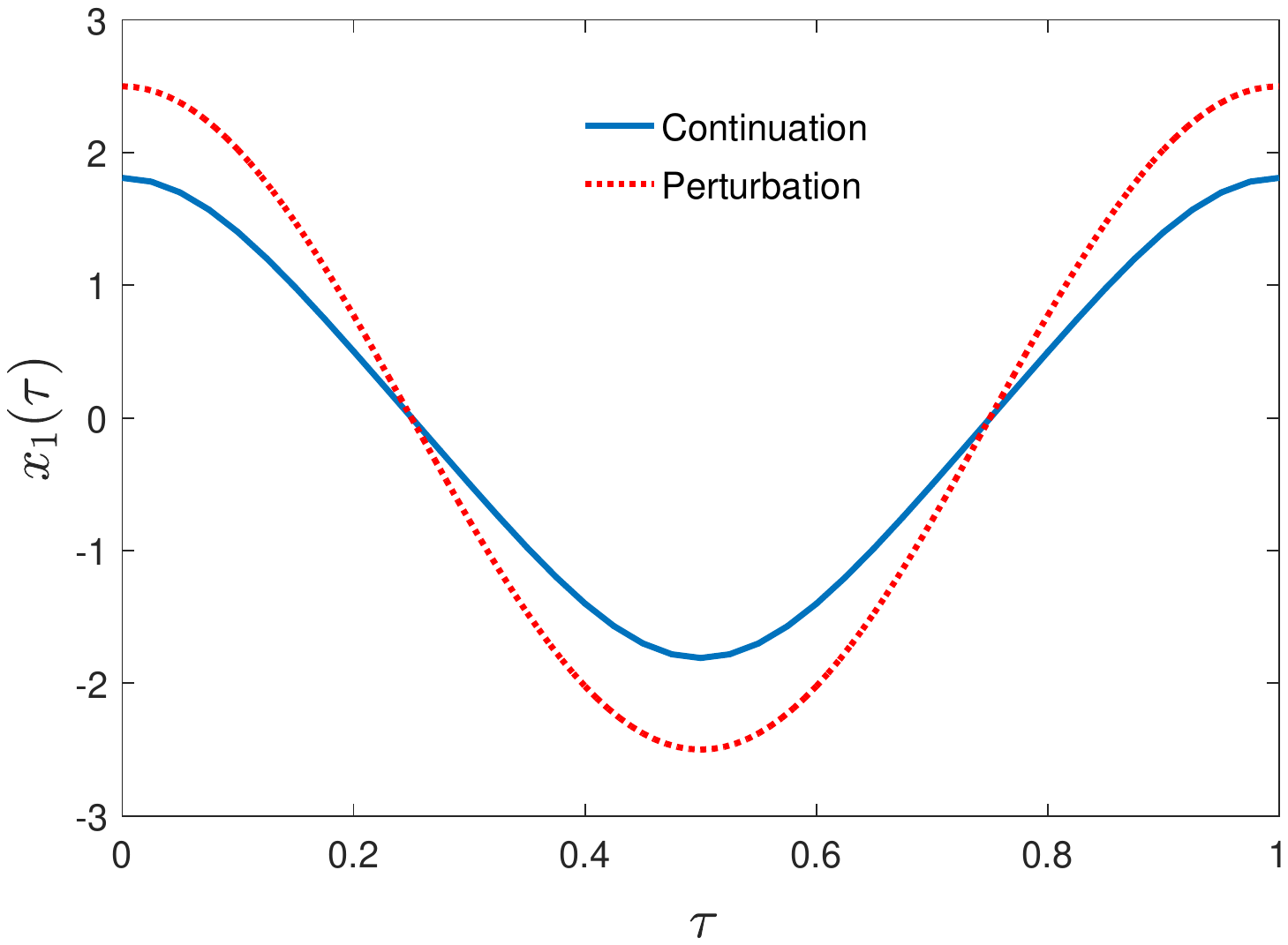}}\\
\subfloat[]{\includegraphics[width=0.9\columnwidth,trim={0 {0.0\textwidth} 0 0.0\textwidth},clip]{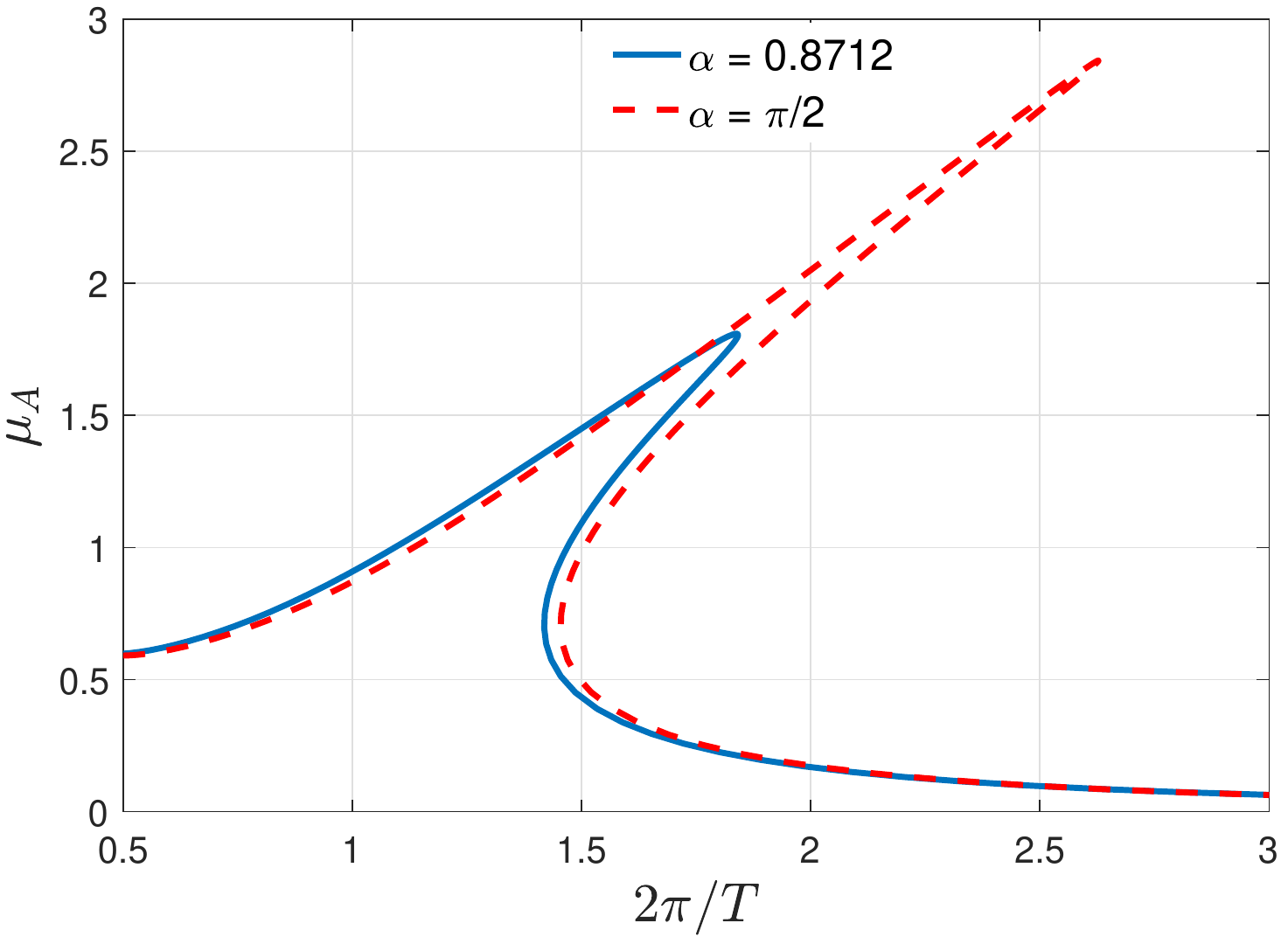}}\\
%\subfloat{\includegraphics[width=0.49\textwidth,trim={0 {0.0\textwidth} 0 0.0\textwidth},clip]{fig_ex2_4c}}
\caption{Optimization of the displacement amplitude along periodic orbits of the harmonically-excited, strongly nonlinear Duffing oscillator with $\zeta=0.05$, $\mu=1$, $a=0.05$, $b=0$, and $\gamma=0.5$ under variations in $\alpha$ and $T$. (a) Comparison of the displacement profile obtained from continuation at the computed optimal delay $\alpha\approx0.8712$ and period $T\approx3.4192$ with the results predicted by perturbation analysis. (b) Frequency-response diagrams for the computed and predicted critical delay values $0.8712$ and $\pi/2$, respectively.}
\label{fig7}
\end{figure}

Further comparisons between the results obtained using the successive continuation approach and those predicted by the perturbation analysis are shown in Figs.~\ref{fig6} and \ref{fig7} for the case when the oscillator is only subjected to displacement feedback, i.e., when $b=0$, with weak ($\mu=0.05$) and strong ($\mu=1$) nonlinearity, respectively. In each case, the perturbation analysis predicts a saddle in the response amplitude for $\alpha=\pi/2\approx 1.5708$, while the computational results are $\alpha\approx 1.4712$ and $0.8712$, respectively. For the case of weak nonlinearity depicted in Fig.~\ref{fig6}, there is still close agreement between the optimal time histories for $x_1(\tau)$, while this is no longer true for the case of strong nonlinearity shown in Fig.~\ref{fig7}. The frequency-response curves shown in the lower panels of Figs.~\ref{fig6} and \ref{fig7} were obtained using numerical continuation for the computed and predicted critical values of $\alpha$. In the case of the weak nonlinearity, we note a weak dependence on the location and magnitude of the peak on the value of the delay, while the differences are stark in the case of the strong nonlinearity. In the latter case, the optimal delay predicted by the perturbation analysis produces a peak amplitude more than $50\%$ larger than that obtained using the numerical method.

\subsection{Quasiperiodic Orbits}

We proceed to consider the problem of optimizing a scalar-valued objective functional on a family of quasiperiodic solutions of \eqref{field}, for which there exists an irrational \emph{rotation number} $\varrho$ and a smooth function $Z:\mathbb{S}\times\mathbb{S}\rightarrow \mathbb{R}^n$ (here, $\mathbb{S}$ denotes the unit circle) such that
\begin{equation}
z(t)=Z\left(\theta_1(t),\theta_2(t)\right),\,\dot{\theta}_1=\varrho\Omega,\,\text{ and }\dot{\theta}_2=\Omega\doteq 2\pi/T
\end{equation}
in terms of the period $T$ of the vector field $f$ in its first argument. Let subscripts $_{\theta_1}$ and $_{\theta_2}$ denote partial derivatives with respect to the corresponding arguments. Substitution into the governing equation then yields the partial differential equation (PDE)
\begin{align}
\label{quasi_eq1}
\varrho\Omega Z_{\theta_{1}}+\Omega Z_{\theta_{2}}= f\left(t,Z,Z\left(\theta_{1}-\varrho\Omega\alpha,\theta_{2}-\Omega\alpha\right),p\right)
\end{align}
on the two-dimensional torus $\mathbb{S}\times\mathbb{S}$. 

We decompose this PDE along its characteristics. To this end, consider the continuous function $V:\mathbb{S}\times[0,1]\rightarrow\mathbb{R}^n$ given by
\begin{equation}
\label{quasi_eq2}
V(\varphi,\tau)\doteq Z(\varphi+2\pi\varrho\tau,2\pi\tau),
\end{equation}
such that $\tau=t/T$, $\theta_1(0)=\varphi$, and without loss of generality $\theta_2(0)=0$. Shifting and wrapping of arguments between and along characteristics will occur several times below. To simplify notation, suppose that $T>\alpha$ and introduce the wrapping operation $W$ for a function $V$ on $\mathbb{S}\times[0,1]$ as
\begin{equation}
    \left[W_a^j V\right](\varphi,\tau)=V(\varphi-2\pi j\rho,\tau-a+j),\,j=-1,0,1
\end{equation}
for $\tau-a+j\in[0,1]$ and all $\varphi\in\mathbb{S}$. It follows by periodicity that 
\begin{equation}
Z(\varphi+2\pi\varrho(\tau-a),2\pi(\tau-a))=\left[W_a^j V\right](\varphi,\tau)
\end{equation}
for $\tau-a+j\in[0,1]$ and all $\varphi\in\mathbb{S}$. Differentiation and use of \eqref{quasi_eq1} then implies that
\begin{align}
\label{quasi_4a}
&V_\tau=Tf\left(T\tau,V,W_{\alpha/T}^1 V,p\right),\,(\varphi,\tau)\in\mathbb{S}\times \left(0,\alpha/T\right),\\
\label{quasi_4b}
&V_\tau=Tf\left(T\tau,V,W_{\alpha/T}^0 V,p\right),\,(\varphi,\tau)\in\mathbb{S}\times\left(\alpha/T,1 \right),
\end{align}
along with the boundary conditions
\begin{equation}
V(\varphi,1)-V(\varphi+2\pi\varrho,0)=0,\,\varphi\in\mathbb{S}.
\label{quasi_4d}
\end{equation}

Equations~\eqref{quasi_4a}--\eqref{quasi_4d} are a family of coupled DDE BVPs in time $\tau$, parametrized by the continuous periodic angle $\varphi$. A family of orbit segments $\mathbb{S}\times[0,1]\ni(\varphi,\tau)\mapsto V(\varphi,\tau)\in\mathbb{R}^n$ solving this family of BVPs then spans the sought quasiperiodic invariant torus. Such a family is unique only up to a shift of its argument $\varphi\in\mathbb{S}$. We isolate a locally unique solution by introducing the integral phase condition
\begin{equation}\label{quasi_4e}
  \int_0^{2\pi}\left(V(\varphi,0)-V^*(\varphi)\right)^\text{T} V^{*}_\varphi(\varphi)\,\mathrm{d}\varphi=0
\end{equation}
in terms of a given continuously-differentiable reference function $V^*:\mathbb{S}\rightarrow\mathbb{R}^n$ that is either fixed throughout the analysis or updated as appropriate. For fixed values of the problem delay $\alpha$, excitation period $T$, and problem parameters $p$, the resultant integro-differential BVP \eqref{quasi_4a}--\eqref{quasi_4e} defining the quasiperiodic response is over-determined (recall that the rotation number $\varrho$ is fixed) such that one has to leave at least one system parameter free to vary to obtain isolated solutions. For example, for fixed $\alpha$, we thus expect to obtain a one-dimensional manifold of quasiperiodic invariant tori under simultaneous variations in $T$ and a single element of $p$.
  
We now apply the construction of the Lagrangian and adjoint equations to this family of DDE BVPs to formulate optimization problems with constraints of the form \eqref{quasi_4a}--\eqref{quasi_4e}, following the procedure from section~\ref{sec_periodic}. We assume that neither the objective functional nor any additional constraints depend on $V$ evaluated for $\tau$ on the interior of the interval $[0,1]$, and that they only depend on $V$ on the boundaries $\tau=0$ and $\tau=1$ through integrals over $\varphi$. In this case, the Lagrange multipliers $\lambda_f$ for the DDE constraint \eqref{quasi_4a} will be continuous on the domain $\mathbb{S}\times[0,1]$ (including periodicity in their first argument $\varphi$). The partial Lagrangian for the constraints \eqref{quasi_4a}--\eqref{quasi_4e} is then given by
\begin{align}
\lefteqn{L_\mathrm{BVP}(V(\cdot,\cdot),\alpha,T,p,\lambda_f(\cdot,\cdot),\lambda_\mathrm{rot}(\cdot),\lambda_\mathrm{ph}) =} \nonumber\\
&\intphi\!\!
%^{2\pi}\!\!\mathrm{d}\varphi
\int\limits_{0}^{\alpha/T}\!\!\!\mathrm{d}\tau
\left[\lambda_{f}^\text{T}\left(V_\tau-Tf_{1}\right)\right]
%\, \nonumber\\ &\quad 
+\intphi\!\!\!
%+\int\limits_{0}^{2\pi}\!\!\mathrm{d}\varphi\!\!\!
\int\limits_{\alpha/T}^{1-\alpha/T}\!\!\!\!\!\!\mathrm{d}\tau\left[\lambda_{f}^\text{T}\left(V_\tau-Tf_{0}\right)\right]
%\,\mathrm{d}\tau\, \mathrm{d}\varphi
\nonumber\\
&+\intphi%\int\limits_{0}^{2\pi}\!\!\mathrm{d}\varphi\!\!\!\!\!
\int\limits_{1-\alpha/T}^{1}\!\!\!\mathrm{d}\tau
\left[\lambda_{f}^\text{T}\left(V_\tau-Tf_{0}\right)\right]\nonumber\\
&+\intphi
%\int\limits_{0}^{2\pi}\!\!\mathrm{d}\varphi
\lambda_{\mathrm{rot}}^\text{T}(\varphi)\left(V\left(\varphi,1\right)-V\left(\varphi+2\pi\varrho,0\right)\right)\,\nonumber\\
&+\lambda_\mathrm{ph}\intphi
%\int\limits_0^{2\pi}\!\!\mathrm{d}\varphi
\left(V(\varphi,0)-V^*(\varphi)\right)^\text{T} V^{*}_\varphi(\varphi),
\label{eq_Lpartial2}
\end{align}
where we abbreviate $\intsphi=\int_0^{2\pi}\mathrm{d}\phi$ and, similarly to section~\ref{sec_periodic}, let $f_j(\varphi,\tau)=f(T\tau, V,W^j_{\alpha/T} V,p)$. The vector-valued functions
$\lambda_{f}(\varphi,\tau)$ and $\lambda_\mathrm{rot}(\varphi)$, and the scalar $\lambda_\mathrm{ph}$ are the Lagrange multipliers associated with the imposition of the differential equations \eqref{quasi_4a} and \eqref{quasi_4b}, boundary conditions \eqref{quasi_4d}, and the integral phase
condition \eqref{quasi_4e}, respectively. Each integrand is assumed to be continuously differentiable on the corresponding interval, and $\lambda_f$ and $\lambda_\mathrm{rot}$ are assumed to be continuous and, hence, periodic in $\varphi$ for all $\tau$. It is again straightforward to show that $\lambda_f$ must be continuous in $\tau$ on $\tau=\alpha/T$ and $\tau=1-\alpha/T$ at a stationary point of the total Lagrangian. In this case, $\lambda_f$ is continuously differentiable in $\tau$ everywhere except at $\tau=1-\alpha/T$, where a slope discontinuity is anticipated from the boundary conditions \eqref{quasi_4d}.

Analogously to Section~\ref{sec_periodic}, consider the notation
\begin{align}
  f_{j,k}(\varphi,\tau)&=\partial_kf(T\tau,V(\varphi,\tau),[W^j_{\alpha/T} V](\varphi,\tau),p),\\
f_{j,q}(\varphi,\tau)&=\frac{\d}{\d\, q}f\left(T\tau,V(\varphi,\tau),[W^j_{\alpha/T} V](\varphi,\tau),p\right),
\end{align}
for $j=0,1$ and $q=\alpha,T$. Then, the contributions to the necessary adjoint conditions for a stationary point of the total Lagrangian are given by
\begin{align}\label{quasi:adj:v1}
&-\lambda_{f,\tau}^\text{T}-T\lambda_{f}^\text{T} f_{1,2}
%\left(T\tau, v,W^1_\alpha v,p\right)\\
-T\left(W^0_{-\alpha/T}\,\lambda_{f}\right)^\text{T} W^0_{-\alpha/T}f_{0,3}
%\left(T\tau+\alpha,W^0_{-\alpha}v,v,p\right)
%g_3 
\end{align}
for variations with respect to $V(\varphi,\tau)$ on $(\varphi,\tau) \in\mathbb{S}\times \left(0,\alpha/T\right)$;
\begin{align}
\label{quasi:adj:v2}
&-\lambda_{f,\tau}^\text{T}-T\lambda_{f}^\text{T} f_{0,2}-T\left(W^0_{-\alpha/T}\,\lambda_{f}\right)^\text{T} W^0_{-\alpha/T} f_{0,3}
%g_3
\end{align}
for variations with respect to $V(\varphi,\tau)$ on $(\varphi,\tau) \in \mathbb{S}\times \left(\alpha/T,1-\alpha/T\right)$;
\begin{align}\label{quasi:adj:v3}
&-\lambda_{f,\tau}^\text{T}-T\lambda_{f}^\text{T} f_{0,2}-T\left(W^{-1}_{-\alpha/T}\,\lambda_{f}\right)^\text{T}
W^{-1}_{-\alpha/T}f_{1,3}
\end{align}
for variations with respect to $V(\varphi,\tau)$ on $(\varphi,\tau) \in \mathbb{S}\times \left(1-\alpha/T,1\right)$;
\begin{equation}
\lambda_{f}^\text{T}\left(\varphi,0\right)+\lambda_\mathrm{rot}^\text{T}\left(\varphi-2\pi\varrho\right)+\lambda_\mathrm{ph} V^{*\top}_\varphi(\varphi)
\label{eq_quasi_tau=0}
\end{equation} 
for variations with respect to  $V(\varphi,0)$ on $\varphi\in\mathbb{S}$;
\begin{equation}
\lambda_{f}^\text{T}\left(\varphi,1\right)+\lambda_\mathrm{rot}^\text{T}\left(\varphi\right)
\label{eq_quasi_tau=1}
\end{equation}
for variations with respect to $V(\varphi,1)$ on $\varphi\in\mathbb{S}$;
\begin{equation}
\label{quasi:adj:alpha}
-\intphi\!\!\int\limits_{0}^{\alpha/T}\!\!\!\mathrm{d}\tau
\left[\lambda_f^\text{T} Tf_{1,\alpha}\right]-
\intphi\int\limits_{\alpha/T}^{1}\!\!\!\mathrm{d}\tau\left[\lambda_f^\text{T} Tf_{0,\alpha}\right]
\end{equation}
for variations with respect to $\alpha$;
\begin{align}
&-\intphi\!\!\int\limits_{0}^{\alpha/T}\!\!\!\mathrm{d}\tau
\left[\lambda_f^\text{T}\left(Tf_{1,T}+f_{1}\right)\right]
%\mathrm{d}\tau\mathrm{d}\varphi\nonumber\\
-\intphi\int\limits_{\alpha/T}^{1}\!\!\!\mathrm{d}\tau
\left[\lambda_f^\text{T}\left(Tf_{0,T}+f_{0}\right)\right]
%\mathrm{d}\tau\mathrm{d}\varphi
\end{align}
for variations with respect to $T$; and
\begin{equation}\label{quasi:adj:p}
-\intphi\int\limits_{0}^{\alpha/T}\!\!\!\mathrm{d}\tau\left[\lambda_f^\text{T} Tf_{1,4}\right]
-\intphi\int\limits_{\alpha/T}^{1}\!\!\!\mathrm{d}\tau\left[\lambda_f^\text{T} Tf_{0,4}\right]
\end{equation}
for variations with respect to $p$.

\subsection{A Hopf unfolding with delay and forcing}

Consider, for example, the problem of finding a local maximum in $\omega$ along a family of quasiperiodic invariant tori of the vector field
\begin{align}
\label{eq_Hopf}
&f\left(t,u,v, p\right)=\left(\begin{array}{c}
-\omega u_{2}
+v_{1}\left(1+r\left(\cos 2\pi t/T-1\right)\right)\\
\omega u_{1}
+v_{2}\left(1+r\left(\cos 2\pi t/T-1\right)\right)
\end{array}\right),
\end{align}
where $r=\sqrt{u_{1}^{2}+u_{2}^{2}}$, $\alpha=1$, and $p=\omega$. Notably, an example in the tutorial for the \textsc{coco} trajectory segment toolbox~\cite{COCO} shows that no such local maximum exists when $\alpha=0$, since then $\omega T=2\pi\varrho$. In the present case, we consider the optimization problem
\begin{equation}
  \label{quasi:obj}
  \mbox{maximize}\quad \mu_\omega=\omega
\end{equation}
subject to the constraints \eqref{quasi_4a}--\eqref{quasi_4e} (the coupled DDEs with boundary conditions and phase condition, depending on $\varphi$). The problem Lagrangian is then given by
\begin{align}
&L(V(\cdot,\cdot),T,p,\mu_\omega,\lambda_f(\cdot,\cdot),\lambda_\mathrm{rot}(\cdot),\lambda_\mathrm{ph},\eta_\omega)\nonumber\\
&\quad =\mu_\omega+\eta_\omega(\omega-\mu_\omega)\nonumber\\
&\qquad+L_\mathrm{BVP}\left(V(\cdot,\cdot),1,T,p,\lambda_f(\cdot,\cdot),\lambda_\mathrm{rot}(\cdot),\lambda_\mathrm{ph}\right),
\end{align}
where $L_\mathrm{BVP}$ is given in \eqref{eq_Lpartial2} and $\eta_\omega$ is the additional Lagrange multiplier. The necessary conditions for an extremum of the total Lagrangian are then given by (i) the original differential equations and boundary conditions \eqref{quasi_4a}--\eqref{quasi_4e}; (ii) the adjoint conditions (excluding \eqref{quasi:adj:alpha}) obtained by appending $\eta_\omega$ to the variation with respect to $p$ \eqref{quasi:adj:p} and setting all the resulting contributions equal to $0$; and (iii) the condition that $\eta_\omega=1$.

As in previous examples, we immediately note that $\lambda_\mathrm{ph}$ must equal $0$ at a stationary point of the Lagrangian, since the objective function is clearly independent of the particular choice of family $(\varphi,\tau)\mapsto V(\varphi,\tau)$ selected by the phase condition. The adjoint boundary conditions \eqref{eq_quasi_tau=0} and \eqref{eq_quasi_tau=1} then imply that
\begin{equation}
\lambda_f^\text{T} (\varphi,1)-\lambda_{f}^\text{T}(\varphi+2\pi\varrho,0)=0.
\end{equation}
Moreover, direct computation using \eqref{eq_Hopf} and the boundary condition \eqref{quasi_4d} shows that
\begin{equation}
f_{0,3}(\varphi,1)-f_{1,3}(\varphi+2\pi\varrho,0)=0.
\end{equation}
It follows from \eqref{quasi:adj:v2} and \eqref{quasi:adj:v3} that
\begin{align}
&\lambda_{f,\tau}^\text{T}(\varphi,1-1/T)_+ - \lambda_{f,\tau}^\text{T}(\varphi,1-1/T)_- =0,
\end{align}
i.e., that $\lambda_f$ is continuously differentiable in $\tau$ on the entire interval $[0,1]$.

We proceed to locate an extremum by applying the successive continuation technique to the set of equations obtained by omitting the trivial adjoint condition that $\eta_\omega=1$. To this end, we approximate  $V(\varphi,\tau)$, $\lambda_f(\varphi,\tau)$, and $\lambda_\mathrm{rot}(\varphi)$ by truncated Fourier series in $\varphi$ with $\tau$-dependent Fourier coefficient functions, as appropriate, approximated by continuous piecewise-polynomial interpolants on the interval $[0,1]$. Although we anticipate that $\lambda_\mathrm{ph}$ will equal $0$ throughout continuation, we keep $\lambda_\mathrm{ph}$ as an unknown and monitor its value during continuation. We first continue along a one-dimensional manifold along which the Lagrange multipliers always equal $0$, and then branch switch at a local maximum of $\mu_\omega$ to a secondary branch along which the family $V$ remains unchanged, while the Lagrange multipliers vary linearly in $\eta_\omega$. The solution to the necessary conditions for a local stationary point is then obtained once $\eta_\omega=1$ along the secondary branch.

\begin{figure}[h]
\centering
\includegraphics[width=0.9\columnwidth]{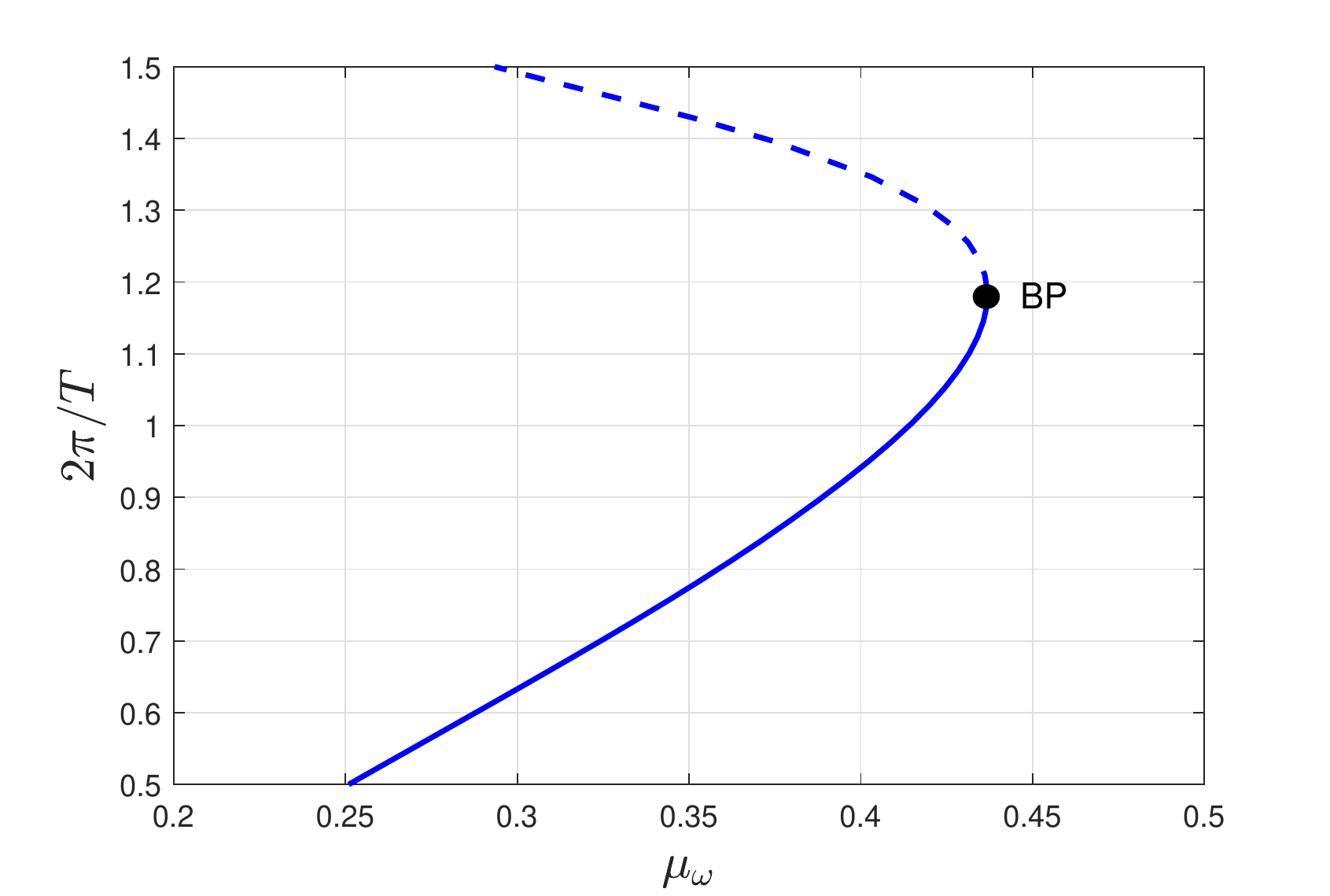}%
\caption{One-dimensional manifold obtained from the first stage of continuation along a family of approximate quasiperiodic invariant tori with vanishing Lagrange multipliers for the case that $\varrho\approx 0.6618$. The local maximum $\mu_\omega\approx0.43685$ when $T\approx 5.3153$ approximately coincides with a branch point (BP). Solid and dashed lines denote dynamically stable and unstable tori, respectively.}%
\label{fig8}
\end{figure}

\begin{figure}[H]
\subfloat[]{\includegraphics[width=.9\columnwidth]{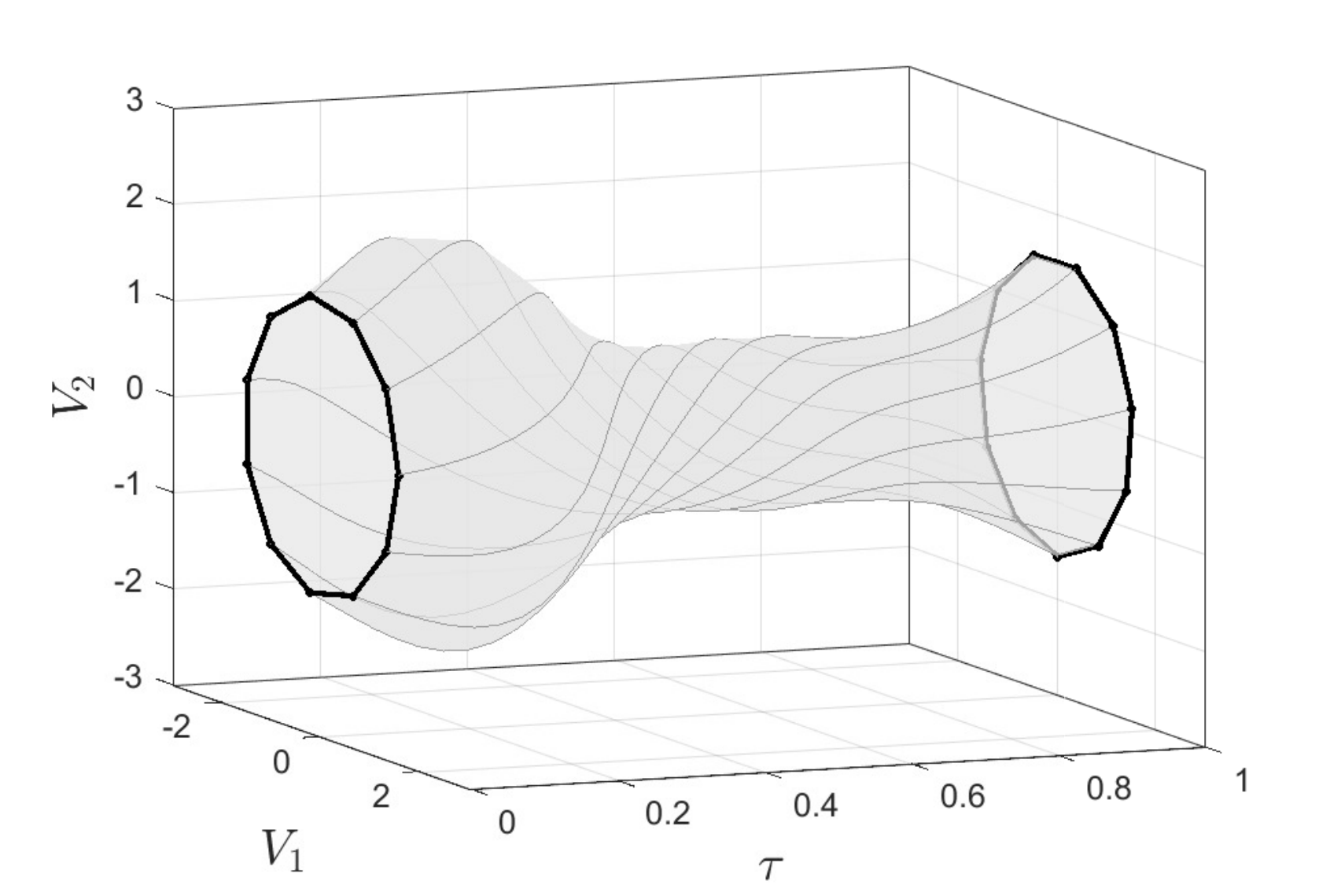}}\\
\subfloat[]{\includegraphics[width=.9\columnwidth]{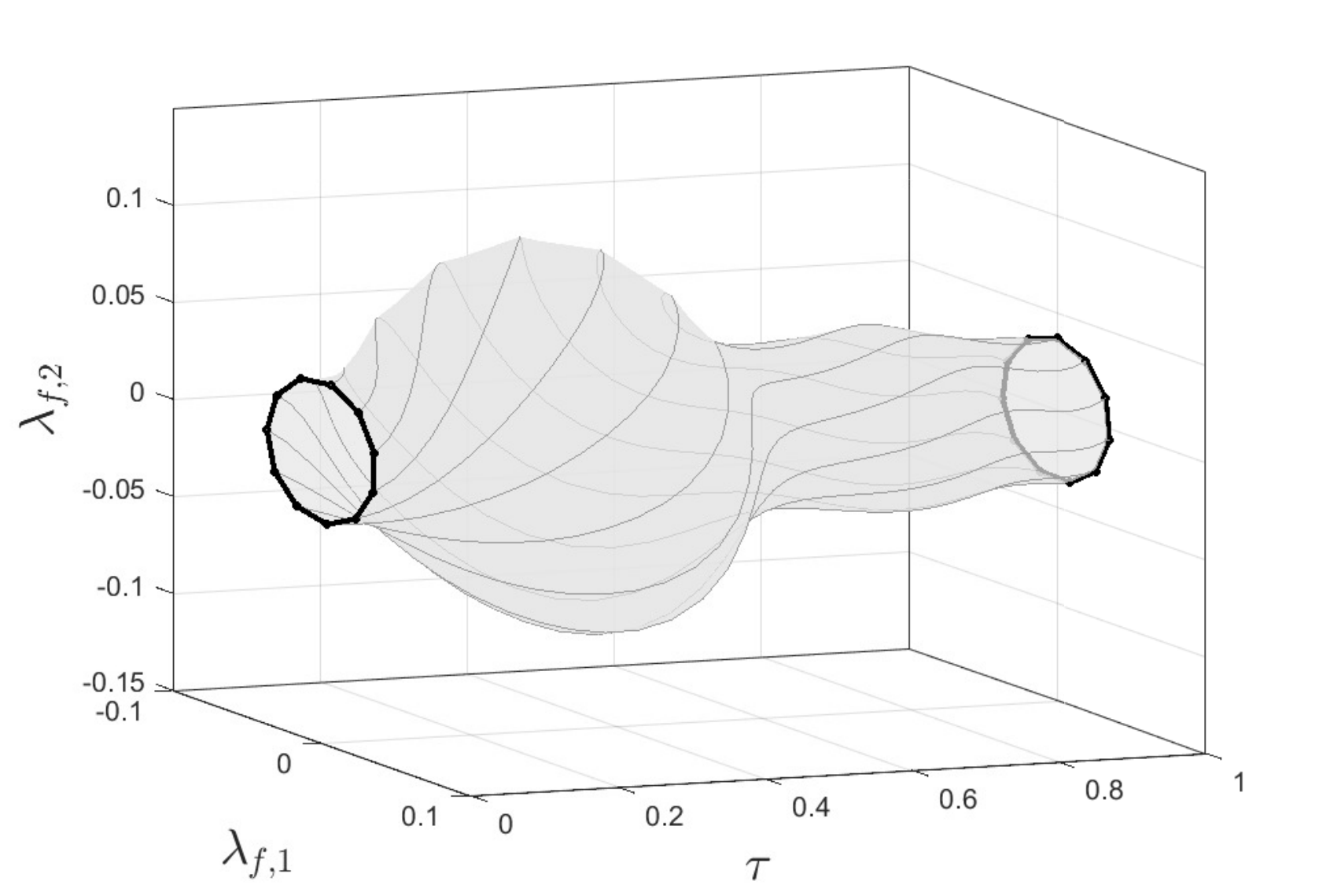}}\\
\subfloat[]{\includegraphics[width=.9\columnwidth]{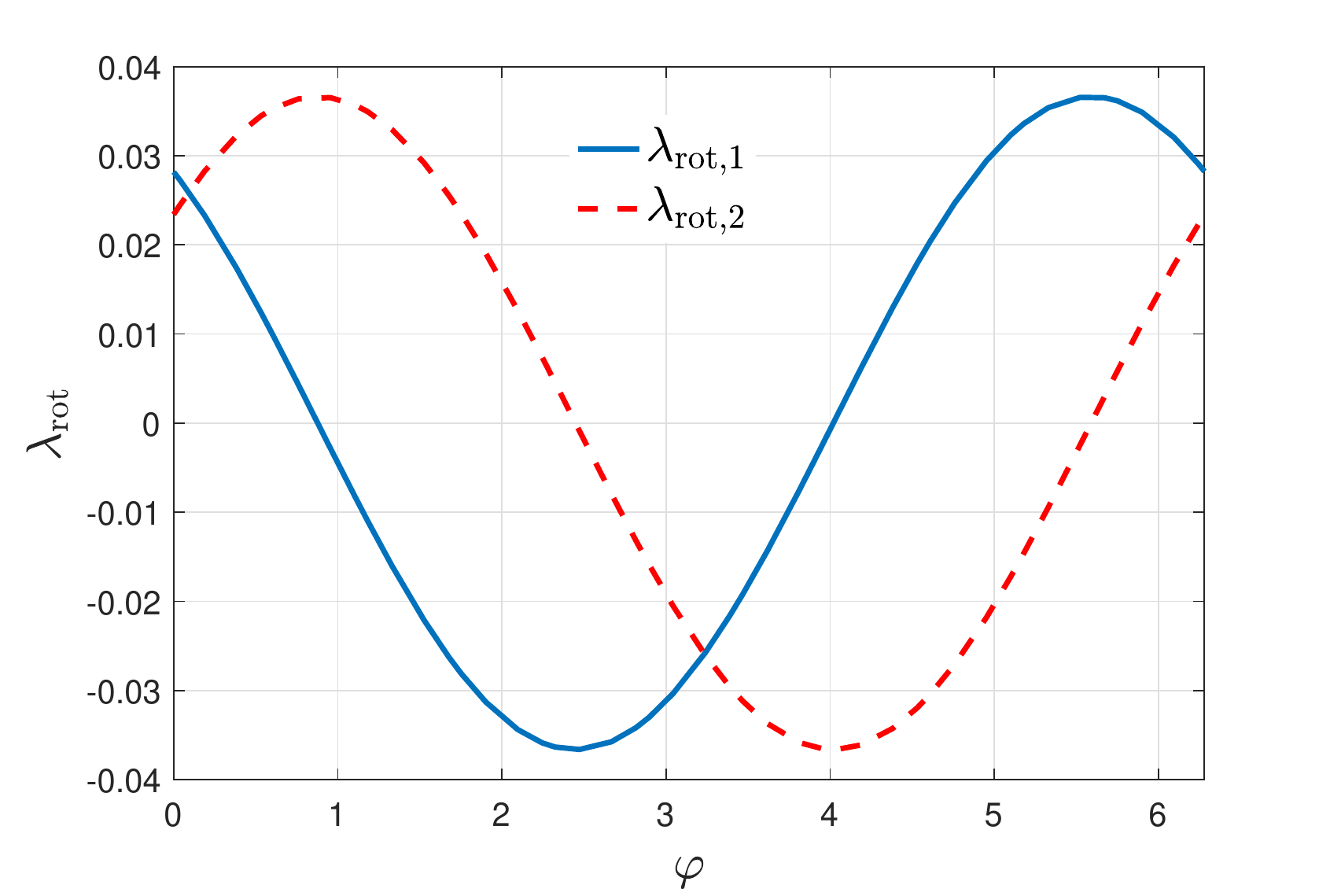}}
\caption{Optimal quasiperiodic invariant torus (a) and corresponding representation of $\lambda_f$ (b) obtained at the terminal point ($\eta_\omega=1$) of the second stage of continuation with $\varrho\approx 0.6618$. Panel (c) shows the components of $\lambda_\mathrm{rot}$ at this point. Solid grey curves in (a) and (b) denote the discretization of $V(\varphi,\tau)$ and $\lambda_f(\varphi,\tau)$ using trajectory segments based at $\varphi=(i-1)/11$, for $i=1,\ldots,11$.}%
\label{fig9}%
\end{figure}

The results of such an analysis using \textsc{coco} is shown in Figs.~\ref{fig8} and \ref{fig9} for the case that $\varrho\approx 0.6618$. Here, dependence on $\varphi$ is approximated using a Fourier series truncated at the fifth harmonic corresponding to $11$ trajectory segments on the torus based at $\varphi=(i-1)/11$, $i=1,\ldots,11$. Each $\tau$ dependent Fourier coefficient is discretized using polynomials of degree $4$ on a uniform mesh with $10$ intervals. The one-dimensional family of quasiperiodic orbits in Fig.~\ref{fig8} along the first manifold with vanishing Lagrange multipliers indicates the existence of a local maximum in $\mu_\omega\approx 0.43685$ for $T\approx5.3153$. Branch switching from the nearby branch point (as before, exact coincidence is lost due to discretization) and continuing until $\eta_\omega=1$ yields the approximate torus and the corresponding Lagrange multipliers $\lambda_f$ and $\lambda_\mathrm{rot}$ shown in Fig.~\ref{fig9}.

As an aside, direct numerical simulation using initial conditions predicted by the continuation analysis suggest that quasiperiodic tori found on the lower half of the one-dimensional family shown in Fig.~\ref{fig8} are stable to sufficiently small perturbations, while the tori found on the upper half are unstable, with a critical loss of stability coincident with the peak value of $\mu_\omega$. 

\section{Conclusions} 
\label{conclusion}

The various examples in previous sections illustrate the successful application to the case with single time delays of the general methodology to optimization along implicitly defined solutions to integro-differential boundary-value problems first proposed by Kernevez and Doedel~\cite{kernevez1987optimization} for ordinary differential equations. Here, the partial Lagrangian approach introduced in \cite{li2017staged} was used to derive adjoint conditions that were linear and homogeneous in the unknown Lagrange multipliers. This allowed a search for local extrema to proceed along a connected sequence of one-dimensional manifolds of solutions to the necessary conditions for such extrema minus the trivial algebraic adjoint conditions on a subset of the Lagrange multipliers: first, along a branch with vanishing Lagrange multipliers, then switching to a branch with linearly varying Lagrange multipliers, and then along additional branches until all the previously omitted trivial algebraic adjoint conditions were satisfied.

In contrast to the case of ordinary differential equations, the presence of time delays introduces potential discontinuities that must be accounted for in any numerical solution strategy. By the properties of differential equations with time-shifted arguments, such discontinuities propagate across time, gaining an order of continuity for each iteration. Here, we have only accounted for zeroth- or first-order discontinuities in the formulation of the governing boundary-value problems. On each segment along which a function was shown to be continuously differentiable, we have approximated such a function by a continuous piecewise-polynomial function of degree $4$ in each mesh interval, ignoring continuity in the first derivative across mesh boundaries or discontinuities of order two or higher within each mesh interval. The piecewise-polynomial approximants have been used to impose a discretization of the governing differential equations at a set of collocation nodes within each interval and to evaluate functions with time-shifted arguments on the same or other intervals. Such a collocation strategy is consistent with the approach in \cite{engelborghs2002stability}, and there compared to an alternative mesh strategy that depends on the delay. We have not undertaken a detailed analysis of the sensitivity of the results to the numerical mesh or polynomial degree. Notably, while we rely in this paper invariably on uniform meshes, it is common to consider adaptive meshes for which the number of intervals and their relative size may change during continuation. We leave such an implementation for future work.

In all the examples, a Lagrange multiplier associated with a phase condition was found to equal $0$ on a local extremum of the corresponding Lagrangian. As stated previously, we nevertheless retained this Lagrange multiplier as an unknown and monitored its value during continuation. Experiments with the number of mesh intervals were used to determine whether this value was effectively $0$ also in the computational analysis. An alternative would have been to eliminate this variable from the set of adjoint equations while simultaneously eliminating one of the adjoint conditions. In a single instance, this may indeed be useful, but when relying on a general-purpose implementation as envisioned in a planned future implementation of \textsc{coco}, it is better to retain the variable and use its numerical value as an indicator of the accuracy of the solution.

There are a number of directions to go in future work. These include consideration of circumstances in which the ratio $\alpha/T$ violates one or several of the inequalities assumed in the previous sections during continuation. Such violations may necessitate a piecewise definition of the Lagrangian across parameter space with different segmentations of the governing differential equations in each region. Problems with multiple delays, as well as problems with state- or time-dependent delays could also be explored as motivated by particular applications.

\begin{acknowledgements}
We would like to thank Mingwu Li for insightful discussions during the course of this study. J.S.\ gratefully acknowledges the financial support of the EPSRC via grants EP/N023544/1 and EP/N014391/1, and from the European Union's Horizon 2020 Research and Innovation Programme under the Marie Sklodowska-Curie grant agreement No 643073.
\end{acknowledgements}

\bibliographystyle{spmpsci}      % mathematics and physical sciences
\bibliography{references}   % name your BibTeX data base

%\begin{thebibliography}{}
%%
%% and use \bibitem to create references. Consult the Instructions
%% for authors for reference list style.
%\end{thebibliography}
\section*{Appendix}
We review the application of the method of multiple scales to the optimal selection of a time delay that results in a minimal peak displacement amplitude in the harmonically-forced response of a Duffing oscillator under delayed displacement and velocity feedback, as discussed in Sect.~\ref{sec_periodic}.

Consider the delay-differential equation
\begin{align}
&\ddot{z}\left(t\right)+2\epsilon\zeta\dot{z}\left(t\right)+z\left(t\right)+\epsilon\mu z^{3}\left(t\right)=2\epsilon a z\left(t-\alpha\right)\nonumber\\
&\qquad \qquad\qquad+2\epsilon b \dot{z}\left(t-\alpha\right)+\epsilon\gamma\,\cos \left((1+\epsilon\sigma) t\right)
\end{align}
for $0<\epsilon\ll1$. We seek an approximate solution of the form
\begin{align}
z\left(t\right)=z_{0}\left(T_{0},T_{1},\ldots\right)+\epsilon z_{1}\left(T_{0},T_{1},\ldots\right)+\cdots,
\end{align}
where $T_i=\epsilon^i t$. To leading order in $\epsilon$,
\begin{equation}
\label{ex2_mms_eq5}
z_{0}\left(T_{0},T_{1},\ldots\right)=A\left(T_{1},\ldots\right)e^{iT_{0}}+ cc,
\end{equation}
where $cc$ denotes complex conjugate terms. Elimination of secular terms at higher orders in $\epsilon$ then yields a set of conditions on the derivatives of the complex amplitude $A$ with respect to its arguments. In particular, if we let
\begin{equation}
A\left(T_{1},\ldots\right)=\frac{1}{2}\rho(T_{1},\ldots)  e^{i\sigma T_1-\varphi(T_{1},\ldots)},
\end{equation}
it follows from the first-order analysis that steady-state oscillations with angular frequency $1+\epsilon\sigma$ result provided that
\begin{align}
\frac{1}{2} \gamma \sin\varphi&=\zeta  \rho+a \rho  \sin\alpha-b \rho  \cos \alpha,\\
\frac{1}{2}\gamma \cos \varphi&=\rho\left(\sigma  +a \cos \alpha+b\sin \alpha -\frac{3 \mu  \rho ^2}{8}\right).
\end{align}
Elimination of $\varphi$ yields the desired, implicit, frequency-amplitude relationship
\begin{align}
\label{eq_freqamp}
& \rho^{2}\left(\sigma+a\cos \alpha+b\sin \alpha-\frac{3\mu\rho^2}{8}\right)^2+\nonumber\\
&\quad\quad\rho^{2}\left(\zeta+a\sin \alpha-b\cos \alpha\right)^2-\frac{\gamma^2}{4}=0,
\end{align}
from which we deduce the maximum value of $\rho$ given by
\begin{equation}
\rho_\mathrm{max}\doteq\frac{\gamma}{2|\zeta+a\sin\alpha-b\cos\alpha|}
\end{equation}
obtained when 
\begin{equation}
\sigma=\frac{3\mu\rho_\mathrm{max}^2}{8}-a\cos\alpha-b\sin\alpha.
\end{equation}
In the special case that $b=-a$, the maximum value of $\rho$ achieves the local minimum $\gamma/2(\zeta+\sqrt{2}a)^2$ for $\alpha=\pi/4$, while for $b=0$, the local minimum $\gamma/2(\zeta+a)^2$ is obtained when $\alpha=\pi/2$.
\end{document}